\documentclass[12pt]{article}% CTP Chernoff paper

	\usepackage[english]{babel}
	\usepackage{blindtext}
\RequirePackage[round,semicolon,authoryear,sort]{natbib}
\RequirePackage[backref=page,colorlinks,citecolor=blue,urlcolor=blue]{hyperref}

\usepackage{graphics}
\usepackage{graphicx}
\usepackage{caption}
\usepackage{subcaption}
\usepackage{verbatim}
\usepackage{pdfpages}

\usepackage{bbm}
\usepackage{parskip}
\usepackage{mathrsfs}
\usepackage{fullpage}
\usepackage{amsfonts,latexsym,eucal,amsmath,amsthm,amssymb,bm,pdfpages}
\usepackage{array}
\usepackage[T1]{fontenc}
\usepackage[ansinew]{inputenc}
\usepackage{enumerate}
\usepackage{booktabs}
\newcommand{\R}{\mathbb{R}}

\newcommand{\bA}{\boldsymbol{A}}
\newcommand{\bB}{\boldsymbol{B}}

\newcommand{\bD}{\boldsymbol{D}}

\newcommand{\bI}{\boldsymbol{I}}

\newcommand{\bL}{\boldsymbol{L}}
\newcommand{\bM}{\boldsymbol{M}}

\newcommand{\bP}{\boldsymbol{P}}
\newcommand{\bQ}{\boldsymbol{Q}}

\newcommand{\bS}{\boldsymbol{S}}

\newcommand{\bU}{\boldsymbol{U}}

\newcommand{\bW}{\boldsymbol{W}}
\newcommand{\bX}{\boldsymbol{X}}
\newcommand{\bY}{\boldsymbol{Y}}

\newcommand{\bDel}{\boldsymbol{\Delta}}
\newcommand{\bN}{\boldsymbol{\mathcal{V}}}
\newcommand{\bPi}{\boldsymbol{\pi}} %is actually lower-case ...
\newcommand{\bSig}{\boldsymbol{\Sigma}}
\newcommand{\bLambda}{\boldsymbol{\Lambda}}

\newcommand{\bOne}{\boldsymbol{1}}
\newcommand{\bZero}{\boldsymbol{0}}
\newcommand{\bTau}{\boldsymbol{\tau}}

\newcommand{\indefO}{\mathbb{O}(d^{+},d^{-})}
\newcommand{\indefI}{\bI^{d^{+}}_{d^{-}}}
\newcommand{\defO}{\mathbb{O}}

\newcommand{\rmOp}{\textrm{op}}

\newcommand{\rmTr}{\textrm{tr}}
\newcommand{\rmDet}{\textrm{det}}
\newcommand{\rmDiag}{\textrm{diag}}
\newcommand{\rmVec}{\textrm{vec}}
\newcommand{\rmRank}{\textrm{rank}}

\newcommand{\rmSBM}{\textrm{SBM}}

\newcommand{\authorER}{\textnormal{Erd\H{o}s--R\'{e}nyi}}

\newcommand{\probO}{O_{\mathbb{P}}}
\newcommand{\probOLog}{O_{\mathbb{P},\log}}
\newcommand{\probTh}{\Theta_{\mathbb{P}}}

\newcommand{\rbox}[1]{\boxed{\text{{\color{red}#1}}}}

\newcommand{\tn}{\textnormal}

\renewcommand{\tilde}{\widetilde}
\renewcommand{\hat}{\widehat}
\renewcommand{\Pr}{\mathbb{P}}

\newcommand{\Ex}{\mathbb{E}}

\newtheorem{theorem}{Theorem}

\newtheorem{lemma}[theorem]{Lemma}

\theoremstyle{definition}
\newtheorem{remark}{Remark}
\newtheorem{definition}[theorem]{Definition}

\makeatletter
\renewcommand{\maketag@@@}[1]{\hbox{\m@th\normalsize\normalfont#1}}%
\makeatother

\graphicspath{{./}{}}%
\graphicspath{{./Figures/}{}}%
\graphicspath{{./INIslides/}{}}%

\newcommand{\blind}{0}

\begin{document}
	
	\def\spacingset#1{\renewcommand{\baselinestretch}%
		{#1}\small\normalsize} \spacingset{1}
	
	\if0\blind
	{
		\title{
			\bf{On spectral embedding performance and\\ elucidating network structure in\\ stochastic block model graphs}
		}
		\author{
			Joshua Cape and Minh Tang and Carey E.~Priebe\\ 
			Department of Applied Mathematics and Statistics\\
			Johns Hopkins University
		}
		\maketitle
	} \fi
	
	%\smallskip
	\begin{abstract}
		Statistical inference on graphs often proceeds via spectral methods involving low-dimensional embeddings of matrix-valued graph representations, such as the graph Laplacian or adjacency matrix. In this paper, we analyze the asymptotic information-theoretic relative performance of Laplacian spectral embedding and adjacency spectral embedding for block assignment recovery in stochastic block model graphs by way of Chernoff information. We investigate the relationship between spectral embedding performance and underlying network structure (e.g.~homogeneity, affinity, core-periphery, (un)balancedness) via a comprehensive treatment of the two-block stochastic block model and the class of $K$-block models exhibiting homogeneous balanced affinity structure. Our findings support the claim that, for a particular notion of sparsity, loosely speaking, ``Laplacian spectral embedding favors relatively sparse graphs, whereas adjacency spectral embedding favors not-too-sparse graphs.'' We also provide evidence in support of the claim that ``adjacency spectral embedding favors core-periphery network structure.''
	\end{abstract}
	%\vfill
	%\thispagestyle{empty}
	%\newpage
	%\clearpage
	%\setcounter{page}{1}
	%\newpage
	\spacingset{1.45} % DON'T change the spacing!
	
	% XXXXXXXXXXXXXXXXXXXXXXXXXXXXXXXXXXXXXXXXXXXXXXXXXXXXXXXXXXXXXXXXXXXXXXXX
	% XXXXXXXXXXXXXXXXXXXXXXXXXXXXXXXXXXXXXXXXXXXXXXXXXXXXXXXXXXXXXXXXXXXXXXXX
	% XXXXXXXXXXXXXXXXXXXXXXXXXXXXXXXXXXXXXXXXXXXXXXXXXXXXXXXXXXXXXXXXXXXXXXXX
	\newpage
	\section{Preface}
	\label{Sec:Preface}
	The stochastic block model (SBM) \citep{holland1983stochastic} is a simple yet ubiquitous network model capable of capturing community structure that has been widely studied via spectral methods in the mathematics, statistics, physics, and engineering communities. Each vertex in an $n$-vertex $K$-block SBM graph belongs to one of the $K$ blocks (communities), and the probability of any two vertices sharing an edge depends exclusively on the vertices' block assignments (memberships).
	
	This paper provides a detailed comparison of two popular spectral embedding procedures by synthesizing recent advances in random graph limit theory. We undertake an extensive investigation of network structure for stochastic block model graphs by considering sub-models exhibiting various functional relationships, symmetries, and geometric properties within the inherent parameter space consisting of block membership probabilities and block edge probabilities. We also provide a collection of figures depicting relative spectral embedding performance as a function of the SBM parameter space for a range of sub-models exhibiting different forms of network structure, specifically homogeneous community structure, affinity structure, core-periphery structure, and (un)balanced block sizes (see Section~\ref{sec:ElucidatingNetworkStructure}).
	
	The rest of this paper is organized as follows.
	\begin{itemize}
		\item Section~\ref{sec:Introduction} introduces the formal setting considered in this paper and contextualizes this work with respect to the existing statistical network analysis literature.
		\item Section~\ref{sec:Preliminaries} establishes notation, presents the generalized random dot product graph model of which the stochastic block model is a special case, defines the adjacency and Laplacian spectral embeddings, presents the corresponding spectral embedding limit theorems, and specifies the notion of sparsity considered in this paper.
		\item Section~\ref{sec:SpectralEmbeddingPerformance} motivates and formulates a measure of large-sample relative spectral embedding performance via Chernoff information.
		\item Section~\ref{sec:ElucidatingNetworkStructure} presents a treatment of the two-block SBM and certain $K$-block SBMs whereby we elucidate the relationship between spectral embedding performance and network model structure.
		\item Section~\ref{sec:Summary} offers further discussion and some concluding remarks.
		\item Section~\ref{sec:Appendix} provides additional details intended to supplement the main body of this paper.
	\end{itemize}
	 	
	\section{Introduction}
	\label{sec:Introduction}
	Formally, we consider the following stochastic block model setting.
	\vspace{.3em}
	\begin{definition}[$K$-block stochastic block model (SBM)]
		\label{def:SBM}
		Let $K \ge 2$ be a positive integer and $\bPi$ be a vector in the interior of the $(K-1)$-dimensional unit simplex in $\R^{K}$. Let $\bB \in (0,1)^{K \times K}$ be a symmetric matrix with distinct rows. We say $(\bA, \bTau) \sim \rmSBM(\bB, \bPi)$ with scaling factor $ 0 < \rho_{n} \le 1$ provided the following conditions hold. Firstly, $\bTau \equiv (\tau_{1}, \dots, \tau_{n})$ where $\tau_{i}$ are independent and identically distributed (i.i.d.)~random variables with $\mathbb{P}[\tau_{i}=k]=\pi_{k}$. Then, $\bA \in \{0,1\}^{n \times n}$ denotes a symmetric (adjacency) matrix such that, conditioned on $\bTau$, for all $i \le j$, the entries $\bA_{ij}$ are independent Bernoulli random variables with $\mathbb{E}[\bA_{ij}]=\rho_{n}\bB_{\tau_{i},\tau_{j}}$. If only $\bA$ is observed, namely if $\bTau$ is integrated out from $(\bA, \bTau)$, then we write $\bA \sim \rmSBM(\bB, \bPi)$.\footnote{The distinct row assumption removes potential redundancy with respect to block connectivity and labeling. Namely, if rows $k$ and $k^{\prime}$ of $\bB^{\prime}$ are identical, then their corresponding blocks are indistinguishable and can without loss of generality be merged to form a reduced block edge probability matrix $\bB$ with corresponding combined block membership probability $\pi_{k}+\pi_{k^{\prime}}$. We also remark that Definition~\ref{def:SBM} implicitly permits vertex self-loops, a choice that we make for mathematical expediency. Whether or not self-loops are disallowed does not alter the asymptotic results and conclusions presented here.}
		 \hfill $\blacktriangle$
	\end{definition}
	The SBM is an example of an inhomogeneous \authorER\ random graph model \citep{bollobas2007phase} and reduces to the classical \authorER\ model \citep{erdos1959random} in the degenerate case when all the entries of $\bB$ are identical.
	This model enjoys an extensive body of literature focused on spectral methods \citep{vonLuxburg2007tutorial} for statistical estimation, inference, and community detection, including
	\cite{fishkind2013consistent,mcsherry2001spectral,lei2015consistency,rohe2011spectral,sussman2012consistent,sarkar2015role}.
	Considerable effort has also been devoted to the information-theoretic and computational investigation of the SBM as a result of interest in the community detection problem; for an overview see \cite{abbe2018survey}.
	Popular variants of the SBM include the mixed-membership stochastic block model \citep{airoldi2008mixed} and the degree-corrected stochastic block model \citep{karrer2011stochastic}.
	
%	In addition, much is known about the fundamental limits for community detection in the SBM for various spectral and non-spectral methods, both for spectral and non-spectral with respect to fundamental limits of community detection, ph For an overview and discussion of recent developments pertaining to information-theoretic and algorithmic recovery limits for community detection in SBM graphs via spectral (and non-spectral) methods, see the survey \citep{abbe--survey--2017}. For a more general treatment of spectral clustering methodologies, see \citep{vonLuxburg2007tutorial} and the references therein.	
%	Indeed, it well-known that for a broad class of stochastic block model regimes, spectral methods are able to recover latent block vertex membership. Various notions of consistency (Rohe et al. 2011). Much is known for various sparsity regimes in which the model parameters change with the vertex size parameter $n$.
	
	Within the statistics literature, substantial attention has been paid to the class of $K$-block SBMs with positive semidefinite block edge probability matrices $\bB$. This is due in part to the extensive study of the \emph{random dot product graph} (RDPG) \emph{model} \citep{nickel2006,young2007random,athreya2018survey}, a latent position random graph model \citep{hoff2002latent} which includes positive semidefinite SBMs as a special case. Notably, it was recently shown that for the random dot product graph model, both Laplacian spectral embedding (LSE; see Definition~\ref{def:ASE_and_LSE}) and adjacency spectral embedding (ASE; see Definition~\ref{def:ASE_and_LSE}) behave approximately as random samples from Gaussian mixture models \citep{athreya2016limit,tang2016limit}. In tandem with these limit results, the concept of Chernoff information \citep{chernoff1952} was employed in \cite{tang2016limit} to demonstrate that neither Laplacian nor adjacency spectral embedding dominates the other for subsequent inference as a spectral embedding method when the underlying inference task is to recover vertices' latent block assignments. In doing so, the results in \cite{tang2016limit} clarify and complete the groundbreaking work in \cite{sarkar2015role} on comparing spectral clusterings for stochastic block model graphs.
	
	In \cite{tang2016limit} the authors leave open the problem of comprehensively investigating Chernoff information as a measure of relative spectral embedding performance for stochastic block model graphs. Moreover, they do not investigate how relative spectral embedding performance corresponds to underlying network model structure. This is understandable, since the positive semidefinite restriction on $\bB$ limits the possible network structure that can be investigated under the random dot product graph model.
	
	More recently, the limit theory in \cite{tang2016limit} was extended in \cite{rubindelanchy2017generalized} to hold for \emph{all} SBMs within the more flexible framework of the \emph{generalized random dot product graph} (GRDPG) \emph{model}. These developments now make it possible to conduct a more comprehensive Chernoff-based analysis, and that is precisely the aim of this paper. We set forth to formulate and analyze a criterion based on Chernoff information for quantifying relative spectral embedding performance which we then further consider in conjunction with underlying network model structure. The investigation carried out in this paper is, to the best of our knowledge, among the first of its kind in the study of statistical network analysis and random graph inference.
	
	This paper focuses on the following two models which have garnered widespread interest (e.g.~see \cite{abbe2018survey} and the references therein).
	\begin{enumerate}
		\item The two-block SBM with
		$\bB = \Bigl[
				\begin{smallmatrix}
					a & b \\
					b & c
			\end{smallmatrix}
			\Bigr]$
		and $\bPi = (\pi_{1}, 1-\pi_{1})$ where $a,b,c,\pi_{1} \in (0,1)$;
		\item The $K \ge 2$ block SBM exhibiting homogeneous balanced affinity structure, i.e.~$\bB_{ij}=a$ for all $i=j$, $\bB_{ij}=b$ for all $i\neq j$, $0<b<a<1$, and $\bPi = (\tfrac{1}{K},\dots,\tfrac{1}{K})\in\R^{K}$.
	\end{enumerate}
	Using the concept of Chernoff information (Section~\ref{sec:SpectralEmbeddingPerformance}), we obtain an information-theoretic summary characteristic
	$\rho^{\star} \equiv \rho^{\star}(\bB,\bPi)$
	such that the cases $\rho^{\star} > 1$, $\rho^{\star} < 1$, and $\rho^{\star} = 1$ correspond to the preference of spectral embedding procedure based on approximate large-sample relative performance, summarized as ASE~$>$~LSE, ASE~$<$~LSE, and ASE~$=$~LSE, respectively. The above models' low-dimensional parameter spaces facilitate visualizing and analyzing the relationship between network structure (i.e.~$\rmSBM(\bB,\bPi)$) and embedding performance (i.e.~$\rho^{\star}(\bB,\bPi)$).
	
	%The boundary values $0$ and $1$ are excluded in order to avoid possible singularities in the limit theory \citep{delanchy--GRDPG--2017,tang2016limit}.
	
	This paper considers the task of performing inference on a single large graph. As such, we interpret the notion of \emph{sparsity} in reference to the magnitudes of probability parameters, namely the magnitudes of the entries of $\bB$. This notion of sparsity corresponds to the interpretation and intuition of a practitioner wanting to do statistics with an observed graph. We shall, with this understanding in mind, subsequently demonstrate that LSE is preferred as an embedding method in relatively sparse regimes, whereas ASE is preferred as an embedding method in not-too-sparse regimes.
	
	By way of contrast, the scaling factor $\rho_{n}$ in Definition~\ref{def:SBM}, which is included for the purpose of general presentation, indexes a sequence of models wherein edge probabilities change with $n$. We take $\rho_{n}$ to be constant in $n$ which by rescaling is equivalent to setting $\rho_{n} \equiv 1$. Limit theorems are known for regimes where $\rho_{n} \rightarrow 0$ as $n \rightarrow\infty$, but these regimes are uninteresting for single graph inference from the perspective of relative spectral embedding performance \citep{tang2016limit}.

	%Throughout this paper, we draw attention to the \emph{latent} aspect of the GRDPG model, a particular \emph{latent position random graph model}. It is this interpretation that has facilitated the limit theory and subsequently the applicability of Chernoff information for the purposes of comparing ASE to LSE. We provide some discussion in this direction, particularly as relates to the underlying geometric considerations of interest in latent space.

	%We will use the terms ``sparse'' and ``dense'' colloquially, as a practitioner uses them -- sparse means relatively few edges, and dense means relatively many edges.  In all cases we consider the number of edges to be $|E| = \Theta(n^2)$.
	%\cite{tang_priebe2016} shows that $|E| = o(n^2)$ implies LSE dominates ASE in our Chernoff information sense.  But: ``Sparsity? Pffft!''
% XXXXXXXXXXXXXXXXXXXXXXXXXXXXXXXXXXXXXXX
% XXXXXXXXXXXXXXXXXXXXXXXXXXXXXXXXXXXXXXX
% XXXXXXXXXXXXXXXXXXXXXXXXXXXXXXXXXXXXXXX
	\section{Preliminaries}
	\label{sec:Preliminaries}
	\subsection{Notation}
	\label{sec:Notation}
	In this paper, all vectors and matrices are real-valued. The symbols $:=$ and $\equiv$ are used to assign definitions and to denote formal equivalence, respectively. Given a symmetric positive definite $n \times n$ matrix $\bM$, let $\langle \cdot,\cdot \rangle_{\bM} : \R^{n}\times\R^{n}\rightarrow\R$ denote the real inner product induced by $\bM$. Similarly, define the induced norm as $\|\cdot\|_{\bM}:=\sqrt{\langle\cdot,\cdot\rangle_{\bM}}$. In particular, given the $n \times n$ identity matrix $\bI$, denote the standard Euclidean inner product and Euclidean norm by $\langle \cdot,\cdot\rangle \equiv \langle \cdot,\cdot \rangle_{\bI}$ and $\|\cdot\|_{2}:=\sqrt{\langle\cdot,\cdot\rangle}$, respectively. 
	Given an underlying matrix, $\rmDet(\cdot)$ and $\rmTr(\cdot)$ denote the matrix determinant and matrix trace operator, respectively. Given a diagonal matrix $\bD := \rmDiag(d_{11},d_{22},\dots,d_{nn})\in\R^{n \times n}$, $|\bD|$ denotes the entrywise absolute value (matrix) of $\bD$.
	
	The vector of all ones in $\R^{n}$ is denoted by $\bOne_{n}$, whereas the zero matrix in $\R^{m \times n}$ is denoted by $\bZero_{m,n}$. We suppress the indices for convenience when the underlying dimensions are understood, writing instead $\bOne$ and $\bZero$.
	
	Let $\mathbb{N}:=\{1,2,3,\dots\}$ denote the set of natural numbers so that for $n\in\mathbb{N}$, $[n]:=\{1,2,\dots,n\}$.
	For integers $d^{+} \ge 1$, $d^{-} \ge 0$, and $d:=d^{+}+d^{-} \ge 1$, let $\indefI:=\bI_{d^{+}}\bigoplus(-\bI_{d^{-}}) \in \R^{d \times d}$ be the direct sum (diagonal) matrix with identity matrices $\bI_{d^{+}}\in\R^{d^{+}\times d^{+}}$ and $\bI_{d^{-}}\in\R^{d^{-}\times d^{-}}$ together with the convention that $\bI^{d^{+}}_{0}\equiv\bI_{d^{+}}$.
	%and $\bI^{0}_{d^{-}} \equiv - \bI_{d^{-}}$.
	For example, $\bI^{1}_{1} \equiv \textnormal{diag}(1,-1)\in\R^{2 \times 2}$.
	
	For integers $n \ge d \ge 1$, the set of all $n \times d$ real matrices with orthonormal columns shall be denoted by $\mathbb{O}_{n,d}$. Let $\mathbb{O}(d^{+},d^{-})$ denote the indefinite orthogonal group with signature $(d^{+},d^{-})$, and let $\mathbb{O}_{d^{+}}\equiv\mathbb{O}_{d^{+},d^{+}}\equiv\mathbb{O}(d^{+},0)$ denote the orthogonal group in $\R^{d^{+}\times d^{+}}$. In particular, $\bM \in \mathbb{O}(d^{+},d^{-})$ has the characterization $\bM^{\top}\indefI\bM = \indefI$. In the case of the orthogonal group, this characterization reduces to the relationship $\bM^{\top} \equiv \bM^{-1}$.
		
	% XXXXXXXXXXXXXXXXXXXXXXXXXX
	\subsection{The generalized random dot product graph model}
	\label{sec:GRDPG}
	A growing corpus has emerged within the statistics literature focused on the development of theory and applications for the \emph{random dot product graph} (RDPG) \emph{model} \citep{nickel2006,young2007random}. This latent position random graph model associates to each vertex in a graph an underlying low-dimensional vector. These vectors may be viewed as encoding structural information or attributes possessed by their corresponding vertices. In turn, the probability of two vertices sharing an edge is specified through the standard Euclidean inner (dot) product of the vertices' latent position vectors. While simple in concept and design, this model has proven successful in real-world applications in the areas of neuroscience and social networks \citep{lyzinski2017hierarchical}. On the theoretical side, the RDPG model enjoys some of the first-ever statistical theory for two-sample hypothesis testing on random graphs, both semiparametric \citep{tang2017semiparametric} and nonparametric \citep{tang2017nonparametric}. For more on the RDPG model, see the survey \cite{athreya2018survey} and the references therein.
	
	More recently, the \emph{generalized random dot product graph} (GRDPG) \emph{model} was introduced as an extension of the RDPG model that includes as special cases the mixed membership stochastic block model as well as \emph{all} (single membership) stochastic block models \citep{rubindelanchy2017generalized}. Effort towards the development of theory for the GRDPG model has already raised new questions and produced new findings related to the geometry of spectral methods, embeddings, and random graph inference. The present paper further contributes to these efforts.
	
	%We now formally introduce the GRDPG model. Our presentation is formulated more generally than originally presented in \citep{delanchy--GRDPG--2017} with slightly modified notation and is in the spirit of \citep{tang2016limit}.
	\begin{definition}[The generalized random dot product graph (GRDPG) model]
		\label{def:GRDPG}
		For integers $d^{+}\ge 1$ and $d^{-} \ge 0$ such that $d:=d^{+}+d^{-} \ge 1$,
		let $F$ be a distribution on a
		%convex\footnote{The GRDPG model definition as presented in \citep{delanchy--GRDPG--2017} includes the assumption that $\mathcal{X}$ is convex for the purpose of encompassing and characterizing the (undirected) mixed membership stochastic block model via latent position space. In general, the GRDPG model remains widely useful and flexible even without the assumption of convexity. The results in this paper do not require the convexity assumption.}
		set $\mathcal{X} \subset \R^{d}$ such that
		$\langle \bI^{d^{+}}_{d^{-}}x,y\rangle \in [0,1]$ for all $x,y\in\mathcal{X}$. We say that $(\bX, \bA) \sim \textrm{GRDPG}(F)$ with signature $(d^{+}, d^{-})$ and scaling factor $0 < \rho_{n} \le 1$ if the following hold. Let $X_{1}, \dots, X_{n} \sim F$ be independent and identically distributed random (latent position) vectors with
		\begin{equation}
			\bX := [X_{1}| \cdots |X_{n}]^{\top}\in\R^{n \times d}
			\text{ and }
			\bP := \rho_{n} \bX\bI_{d^{-}}^{d^{+}}\bX^{\top}\in[0,1]^{n \times n}.
		\end{equation}
		For each $i \le j$, the entries $\bA_{ij}$ of the symmetric adjacency matrix $\bA\in\{0,1\}^{n \times n}$ are then generated in a conditionally independent fashion given the latent positions, namely
		\begin{equation}
			\label{eq:AijConditionalProbs}
			\{\bA_{ij}|X_{i},X_{j}\} \sim \textrm{Bernoulli}(\rho_{n}\langle \bI_{d^{-}}^{d^{+}}X_{i},X_{j}\rangle).
		\end{equation}
		In this setting, the conditional probability $\mathbb{P}[\bA|\bX]$ can be computed explicitly as a product of Bernoulli probabilities. \hfill $\blacktriangle$
	\end{definition}
	To reiterate, we consider the regime $\rho_{n} \equiv 1$ and therefore suppress dependencies on $\rho_{n}$ later in the text. When no confusion can arise, we also use adorned versions of the symbol $\rho$ to denote Chernoff-related quantities unrelated to $\rho_{n}$ in a manner consistent with the notation in \cite{tang2016limit} (see Section~\ref{sec:SpectralEmbeddingPerformance}).
	
	When $d^{-}=0$, the GRDPG model reduces to the RDPG model. When the distribution $F$ is a discrete distribution on a finite collection of vectors in $\R^{d}$, then the GRDPG model coincides with the SBM, in which case the $n \times n$ edge probability matrix $\bP$ arises as an appropriate dilation of the $K \times K$ block edge probability matrix $\bB$. Given any valid $\bB \in (0,1)^{K\times K}$ as in Definition~\ref{def:SBM}, there exist integers $d^{+}, d^{-}$, and a matrix $\bX\in\R^{K \times K}$ such that $\bB$ has the (not necessarily unique) factorization $\bB \equiv \bX \indefI \bX^{\top}$, which follows since the spectral decomposition of $\bB$ can be written as $\bB \equiv \bU_{\bB}\bLambda\bU_{\bB}^{\top} = (\bU_{\bB}|\bLambda|^{1/2})\indefI(\bU_{\bB}|\bLambda|^{1/2})^{\top}$. This demonstrates the ability of the GRDPG framework in Definition~\ref{def:GRDPG} to model all possible stochastic block models formulated in Definition~\ref{def:SBM}.
	
	\begin{remark}[Non-identifiability in the GRDPG model]
		\label{rem:GRDPG_identifiability}
		The GRDPG model possess two intrinsic sources of non-identifiability, summarized as
		``uniqueness up to indefinite orthogonal transformations'' and
		``uniqueness up to artificial dimension blow-up''.
		More precisely, for $(\bX, \bA) \sim \textrm{GRDPG}(F)$ with signature $(d^{+},d^{-})$, the following considerations must be taken into account.
		\begin{enumerate}
			\item For any $\bQ \in \mathbb{O}(d^{+},d^{-})$, $(\bX, \bA) \overset{\textrm{d}}{=}  (\bY, \bB)$ whenever $(\bY, \bB) \sim \textrm{GRDPG}(F \circ \bQ)$, where $F \circ \bQ$ denotes the distribution of the latent position vector $Y = \bQ X$ and $\overset{\textrm{d}}{=}$ denotes equality in distribution. This source of non-identifiability cannot be mitigated. See Eq.~(\ref{eq:AijConditionalProbs}).
			\item There exists a distribution $F^{\prime}$ on $\R^{d^{\prime}}$ for some $d^{\prime} > d$ such that $(\bX, \bA) \overset{\textrm{d}}{=}  (\bY, \bB)$ where $(\bY, \bB) \sim \textrm{GRDPG}(F^{\prime})$. This source of non-identifiability can be avoided by assuming, as we do in this paper, that $F$ is non-degenerate in the sense that for $X_{1} \sim F$, the second moment matrix $\mathbb{E}[X_{1} X_{1}^{\top}]\in\R^{d \times d}$ is full rank.
		\end{enumerate}
	\end{remark}
	\begin{definition}[Adjacency and Laplacian spectral embeddings]
		\label{def:ASE_and_LSE}
		Let $\bA \in \{0,1\}^{n \times n}$ be a symmetric adjacency matrix with eigendecomposition $\bA \equiv \sum_{i=1}^{n}\lambda_{i}u_{i}u_{i}^{\top}$ and with ordered eigenvalues $|\lambda_{1}| \ge |\lambda_{2}| \ge \cdots \ge |\lambda_{n}|$ corresponding to orthonormal eigenvectors $u_{1}, u_{2}, \dots, u_{n}$. Given a positive integer $d$ such that $d \le n$, let $\bS_{\bA}:=
		\rmDiag(\lambda_{1},\dots,\lambda_{d})\in\R^{d \times d}$ and $\bU_{\bA}:=[u_{1}|\dots|u_{d}]\in\mathbb{O}_{n,d}$. The \emph{adjacency spectral embedding} (ASE) of $\bA$ into $\R^{d}$ is then defined to be the $n \times d$ matrix $\hat{\bX}:=\bU_{\bA}|\bS_{\bA}|^{1/2}$. The matrix $\hat{\bX}$ serves as a consistent estimator
		for $\bX$ up to indefinite orthogonal transformation as $n \rightarrow \infty$.
		
		Along similar lines, define the normalized Laplacian of $\bA$ as
		\begin{equation}
		\label{eq:Laplacian}
			\mathcal{L}(\bA):=(\rmDiag(\bA\bOne_{n}))^{-1/2}\bA(\rmDiag(\bA\bOne_{n}))^{-1/2} \in \R^{n \times n}
		\end{equation}
		whose eigendecomposition is given by $\mathcal{L}(\bA)\equiv\sum_{i=1}^{n}\tilde{\lambda}_{i}\tilde{u}_{i}\tilde{u}_{i}^{\top}$ with ordered eigenvalues $|\tilde{\lambda}_{1}| \ge |\tilde{\lambda}_{2}| \ge \dots \ge |\tilde{\lambda}_{n}|$ corresponding to orthonormal eigenvectors $\tilde{u}_{1}, \tilde{u}_{2}, \dots, \tilde{u}_{n}$. Given a positive integer $d$ such that $d \le n$, let $\tilde{\bS}_{\bA}:=
		\rmDiag(\tilde{\lambda}_{1},\dots,\tilde{\lambda}_{d})\in\R^{d \times d}$ and let $\tilde{\bU}_{\bA}:=[\tilde{u}_{1}|\dots|\tilde{u}_{d}]\in\mathbb{O}_{n,d}$. The \emph{Laplacian spectral embedding} (LSE) of $\bA$ into $\R^{d}$ is then defined to be the $n \times d$ matrix $\breve{\bX}:=\tilde{\bU}_{\bA}|\tilde{\bS}_{\bA}|^{1/2}$. The matrix $\breve{\bX}$ serves as a consistent estimator for the matrix  $(\rmDiag(\bX\indefI\bX^{\top}\bOne_{n}))^{-1/2}\bX$ up to indefinite orthogonal transformation as $n \rightarrow \infty$.
		\hfill $\blacktriangle$
	\end{definition}
	
	\begin{remark}[Consistent estimation and parametrization involving latent positions]
		\label{rem:GRDPG_consistency}
		The matrices $\bX$ and $(\rmDiag(\bX\indefI\bX^{\top}\bOne_{n}))^{-1/2}\bX$, which are one-to-one invertible transformations of each other, may be viewed as providing different parametrizations of GRDPG graphs. As such, comparing $\hat{\bX}$ and $\breve{\bX}$ as estimators is non-trivial. In order to carry out such a comparison, we subsequently adopt an information-theoretic approach in which we consider a particular choice of $f$-divergence which is both analytically tractable and statistically interpretable in the current setting.
	\end{remark}
	
	For the subsequent purposes of the present work, Theorems~\ref{thrm:GRDPG_CLT_ASE}~and~\ref{thrm:GRDPG_CLT_LSE} (below) state slightly weaker formulations of the corresponding limit theorems obtained in \cite{rubindelanchy2017generalized} for adjacency and Laplacian spectral embedding.
	
	\begin{theorem}[ASE limit theorem for GRDPG, adapted from \cite{rubindelanchy2017generalized}]
		\label{thrm:GRDPG_CLT_ASE}
		Assume the $d$-dimensional GRDPG setting in Definition~\ref{def:GRDPG} with $\rho_{n}\equiv 1$. Let $\hat{\bX}$ be the adjacency spectral embedding into $\R^{d}$ with $i$-th row denoted by $\hat{X}_{i}$. Let $\Phi(\cdot,\bSig)$ denote the cumulative distribution function of the centered multivariate normal distribution in $\R^{d}$ with covariance matrix $\bSig$. Then, with respect to the adjacency spectral embedding, there exists a sequence of matrices $\bQ \equiv \bQ_{n} \in \indefO$ such that, for any $z\in\R^{d}$,
		\begin{equation}
			\Pr\left[\sqrt{n}\left(\bQ\hat{X}_{i}-X_{i}\right) \le z\right]
			\rightarrow \int_{\mathcal{X}}\Phi(z,\bSig(x))dF(x)
		\end{equation}
		as $n \rightarrow \infty$, where for $X_{1} \sim F$,
		\begin{equation*}
			\bSig(x) := \indefI\bDel^{-1}\Ex\left[g(x,X_{1})X_{1}X_{1}^{\top} \right]\bDel^{-1}\indefI,
		\end{equation*}
		with $\bDel := \Ex[X_{1}X_{1}^{\top}]$ and $g(x,X_{1}):=\langle \indefI x, X_{1} \rangle(1-\langle \indefI x, X_{1} \rangle)$.
	\end{theorem}
	\begin{theorem}[LSE limit theorem for GRDPG, adapted from \cite{rubindelanchy2017generalized}]
		\label{thrm:GRDPG_CLT_LSE}
		Assume the $d$-dimensional GRDPG setting in Definition~\ref{def:GRDPG} with $\rho_{n}\equiv 1$. Let $\breve{\bX}$ be the Laplacian spectral embedding into $\R^{d}$ with $i$-th row denoted by $\breve{X}_{i}$. Let $\Phi(\cdot,\bSig)$ denote the cumulative distribution function of the centered multivariate normal distribution in $\R^{d}$ with covariance matrix $\bSig$. Then, with respect to the Laplacian spectral embedding, there exists a sequence of matrices $\tilde{\bQ} \equiv \tilde{\bQ}_{n} \in \indefO$ such that, for any $z\in\R^{d}$,
		\begin{equation}
			\Pr\left[n\left(\tilde{\bQ}\breve{X}_{i}-\tfrac{X_{i}}{\sqrt{\sum_{j}\langle\indefI X_{i},X_{j}\rangle}}\right) \le z \right]
			\rightarrow \int_{\mathcal{X}}\Phi(z,\tilde{\bSig}(x))dF(x)
		\end{equation}
		as $n \rightarrow \infty$, where for $X_{1} \sim F$ and $\mu := \Ex[X_{1}]$,
		\begin{equation*}
			\tilde{\bSig}(x) := \indefI\tilde{\bDel}^{-1}\Ex\left[\tilde{g}(x,X_{1})\left(\tfrac{X_{1}}{\langle\indefI \mu, X_{1}\rangle} - \tfrac{\tilde{\bDel}\indefI x}{2\langle\indefI\mu, x \rangle}\right)\left(\tfrac{X_{1}}{\langle\indefI \mu, X_{1}\rangle} - \tfrac{\tilde{\bDel}\indefI x}{2\langle\indefI\mu, x \rangle}\right)^{\top}\right]\tilde{\bDel}^{-1}\indefI,
		\end{equation*}
		with $\tilde{\bDel}:=\Ex\left[\langle\indefI\mu, X_{1}\rangle^{-1} X_{1}X_{1}^{\top}\right]$ and $\tilde{g}(x,X_{1}):=\left[\langle\indefI\mu,x\rangle^{-1}\langle\indefI x, X_{1}\rangle(1-\langle\indefI x, X_{1}\rangle)\right]$.
	\end{theorem}
% XXXXXXXXXXXXXXXXXXXXXXXXXXXXXXXXXXXXXXXXXXXXX
% XXXXXXXXXXXXXXXXXXXXXXXXXXXXXXXXXXXXXXXXXXXXX
% XXXXXXXXXXXXXXXXXXXXXXXXXXXXXXXXXXXXXXXXXXXXX
	\section{Spectral embedding performance}
	\label{sec:SpectralEmbeddingPerformance}
	We desire to compare the large-$n$ sample relative performance of adjacency and Laplacian spectral embedding for subsequent inference, where the subsequent inference task is naturally taken to be the problem of recovering latent block assignments. Here, measuring spectral embedding performance will correspond to estimating the large-sample optimal error rate for recovering the underlying block assignments following each of the spectral embeddings. Towards this end, we now introduce Chernoff information and Chernoff divergence as appropriate information-theoretic quantities.
	
	Given independent and identically distributed random vectors $Y_{i}$ arising from one of two absolutely continuous multivariate distributions $F_{1}$ and $F_{2}$ on $\Omega = \R^{d}$ with density functions $f_{1}$ and $f_{2}$, respectively, we are interested in testing the simple null hypothesis $\mathbb{H}_{0}: F=F_{1}$ against the simple alternative hypothesis $\mathbb{H}_{\textrm{A}}:F=F_{2}$. In this framework, a statistical test $T$ can be viewed as a sequence of mappings $T_{m}:\Omega^{m}\rightarrow\{1,2\}$ indexed according to sample size $m$ such that $T_{m}$ returns the value two when $\mathbb{H}_{0}$ is rejected in favor of $\mathbb{H}_{\textrm{A}}$ and correspondingly returns the value one when $\mathbb{H}_{0}$ is favored. For each $m$, the corresponding significance level and type-II error are denoted by $\alpha_{m}$ and $\beta_{m}$, respectively.
	
	Assume that the prior probability of $\mathbb{H}_{0}$ being true is given by $\pi\in(0,1)$. For a given $\alpha_{m}^{\star}\in(0,1)$, let $\beta_{m}^{\star}\equiv\beta_{m}^{\star}(\alpha_{m}^{\star})$ denote the type-II error associated with the corresponding likelihood ratio test when the type-I error is at most $\alpha_{m}^{\star}$. Then, the \emph{Bayes risk} in deciding between $\mathbb{H}_{0}$ and $\mathbb{H}_{\textrm{A}}$ given $m$ independent random vectors $Y_{1}, Y_{2},\dots, Y_{m}$ is given by
	\begin{equation}
		\label{eq:Bayes risk}
		\underset{\alpha_{m}^{\star}\in(0,1)}{\textrm{inf}}{\pi\alpha_{m}^{\star}+(1-\pi)\beta_{m}^{\star}.}
	\end{equation}
	The Bayes risk is intrinsically related to \emph{Chernoff information} \citep{chernoff1952,chernoff1956}, $C(F_{1}, F_{2})$, namely
	\begin{equation}
		\label{eq:limitErrorChernoff}
		\underset{m\rightarrow\infty}{\textrm{lim}}\tfrac{1}{m}
		\left[\underset{\alpha_{m}^{\star}\in(0,1)}{\textrm{inf}}{\log(\pi\alpha_{m}^{\star}+(1-\pi)\beta_{m}^{\star})}\right]
			= -C(F_{1}, F_{2}),
	\end{equation}
	where
	\begin{align*}
		C(F_{1}, F_{2})
		&:=  -\log\left[\underset{t\in(0,1)}{\textrm{inf}} \int_{\R^{d}}f_{1}^{t}(\boldmath{x})f_{2}^{1-t}(\boldmath{x})d\boldmath{x}\right]
		=  \underset{t\in(0,1)}{\textrm{sup}}\left[-\log \int_{\R^{d}}f_{1}^{t}(\boldmath{x})f_{2}^{1-t}(\boldmath{x})d\boldmath{x}\right].
	\end{align*}
	In words, the Chernoff information between $F_{1}$ and $F_{2}$ is the exponential rate at which the Bayes risk decreases as $m \rightarrow\infty$. Note that the Chernoff information is independent of the prior probability $\pi$. A version of Eq.~(\ref{eq:limitErrorChernoff}) also holds when considering $K \ge 3$ hypothesis with distributions $F_{1}, F_{2}, \dots, F_{K}$, thereby introducing the quantity $\underset{k \neq l}{\textrm{min}} \ C(F_{k}, F_{l})$ (see for example \cite{tang2016limit}).
	
	Chernoff information can be expressed in terms of the \emph{Chernoff divergence} between distributions $F_{1}$ and $F_{2}$, defined for $t\in(0,1)$ as
	\begin{equation}
		C_{t}(F_{1}, F_{2}) = -\log \int_{\R^{d}}f_{1}^{t}(x)f_{2}^{1-t}(x)dx,
	\end{equation}	
	which yields the relation
	\begin{equation}
		C(F_{1}, F_{2}) = \underset{t\in(0,1)}{\textrm{sup}}C_{t}(F_{1}, F_{2}).
	\end{equation}	
	The Chernoff divergence is an example of an $f$-divergence and as such satisfies the data processing lemma \citep{liese2006divergences} and is invariant with respect to invertible transformations \citep{devroye2013probabilistic}. One could instead use another $f$-divergence for the purpose of comparing the two embedding methods, such as the Kullback-Liebler divergence. Our choice is motivated by the aforementioned relationship with Bayes risk in Eq.~(\ref{eq:limitErrorChernoff}).
	
	In this paper we explicitly consider multivariate normal distributions as a consequence of Theorems~\ref{thrm:GRDPG_CLT_ASE}~and~\ref{thrm:GRDPG_CLT_LSE} when conditioning on the individual underlying latent positions for stochastic block model graphs. In particular, given $F_{1}=\mathcal{N}(\mu_{1},\bSig_{1})$, $F_{2}=\mathcal{N}(\mu_{2},\bSig_{2})$, and $t\in(0,1)$, then for $\bSig_{t}:=t\bSig_{1}+(1-t)\bSig_{2}$, the Chernoff information between $F_{1}$ and $F_{2}$ is given by
	\begin{align*}
		C(F_{1}, F_{2})
		&= \underset{t\in(0,1)}{\textrm{sup}}\left[\tfrac{t(1-t)}{2}(\mu_{2}-\mu_{1})^{\top}\bSig_{t}^{-1}(\mu_{2}-\mu_{1}) + \tfrac{1}{2}\log\left(\tfrac{\rmDet(\bSig_{t})}{\rmDet(\bSig_{1})^{t}\rmDet(\bSig_{2})^{1-t}}\right)\right]\\
		&= \underset{t\in(0,1)}{\textrm{sup}}\left[\tfrac{t(1-t)}{2}\|\mu_{2}-\mu_{1}\|_{\bSig_{t}^{-1}}^{2} + \tfrac{1}{2}\log\left(\tfrac{\rmDet(\bSig_{t})}{\rmDet(\bSig_{1})^{t}\rmDet(\bSig_{2})^{1-t}}\right)\right].
	\end{align*}
	%The above discussion together with Theorems~\ref{thrm:GRDPG_CLT_ASE}--\ref{thrm:GRDPG_CLT_LSE} informs the use of the minimum of pairwise Chernoff informations between appropriate multivariate normal distributions for measuring the large-sample optimal error rate for spectral clustering using adjacency or Laplacian spectral embedding. In particular,
	Let $\bB\in(0,1)^{K \times K}$ and $\bPi$ denote the matrix of block edge probabilities and the vector of block assignment probabilities for a $K$-block stochastic block model as before. This corresponds to a special case of the GRDPG model with signature $(d^{+},d^{-})$, $d^{+}+d^{-}=\rmRank(\bB)$, and latent positions $\nu_{k} \in \R^{\rmRank(\bB)}$. For an $n$-vertex SBM graph with parameters $(\bB, \bPi)$, the large-sample optimal error rate for recovering block assignments when performing adjacency spectral embedding can be characterized by the quantity $\rho_{\textrm{A}}\equiv\rho_{\textrm{A}}(\bB, \bPi, n)$ defined by
	\begin{equation}
	\label{eq:rhoA}
		\rho_{\tn{A}} :=
		\underset{k \neq l}{\tn{min}}\underset{t\in(0,1)}{\tn{sup}}
		\left[
		\tfrac{n t(1-t)}{2}\|\nu_{k}-\nu_{l}\|_{\bSig_{kl}^{-1}(t)}^{2} + \tfrac{1}{2}\log\left(\tfrac{\rmDet(\bSig_{kl}(t))}{\rmDet(\bSig_{k})^{t}\rmDet(\bSig_{l})^{1-t}}\right)
		\right],
	\end{equation}
	where $\bSig_{kl}(t):=t\bSig_{k}+(1-t)\bSig_{l}$ for $t \in (0,1)$.
	
	Similarly, for Laplacian spectral embedding, $\rho_{\tn{L}}\equiv\rho_{\tn{L}}(\bB, \bPi, n)$, one has
	\begin{equation}
	\label{eq:rhoL}
		\rho_{\tn{L}} :=
		\underset{k \neq l}{\tn{min}}\underset{t\in(0,1)}{\tn{sup}}
		\left[
		\tfrac{n t(1-t)}{2}\|\tilde{\nu}_{k}-\tilde{\nu}_{l}\|_{\tilde{\bSig}_{kl}^{-1}(t)}^{2} + \tfrac{1}{2}\log\left(\tfrac{\rmDet(\tilde{\bSig}_{kl}(t))}{\rmDet(\tilde{\bSig}_{k})^{t}\rmDet(\tilde{\bSig}_{l})^{1-t}}\right)
		\right],
	\end{equation}
	where $\tilde{\bSig}_{kl}(t):=t\tilde{\bSig}_{k}+(1-t)\tilde{\bSig}_{l}$ and
	$\tilde{\nu}_{k} := \nu_{k}/(\sum_{k^{\prime}}\pi_{k^{\prime}}\langle \indefI\nu_{k^{\prime}},\nu_{k}\rangle)^{1/2}$.
	
	The factor $n$ in Eqs.~(\ref{eq:rhoA}--\ref{eq:rhoL}) arises from the implicit consideration of the appropriate (non-singular) theoretical sample covariance matrices. To assist in the comparison and interpretation of the quantities $\rho_{\textrm{A}}$ and $\rho_{\textrm{L}}$, we assume throughout this paper that $n_{k}=n\pi_{k}$ for $\tilde{\nu}_{k}$. The logarithmic terms in Eqs.~(\ref{eq:rhoA}--\ref{eq:rhoL}) as well as the deviations of each term $n_{k}$ from $n\pi_{k}$ are negligible for large $n$, collectively motivating the following large-sample measure of relative performance, $\rho^{\star}$, where
	\begin{equation}
		\label{eq:rhoStar}
		\frac{\rho_{\tn{A}}}{\rho_{\tn{L}}}
		\equiv 
		\frac{\rho_{\tn{A}}(n)}{\rho_{\tn{L}}(n)}
		\rightarrow
		\rho^{\star}
		\equiv \frac{\rho_{\tn{A}}^{\star}}{\rho_{\tn{L}}^{\star}}
		:= \frac{\underset{k \neq l}{\tn{min}}\underset{t\in(0,1)}{\tn{sup}}
			\left[
			t(1-t)\|\nu_{k}-\nu_{l}\|_{\bSig_{kl}^{-1}(t)}^{2}\right]}
		{\underset{k \neq l}{\tn{min}}\underset{t\in(0,1)}{\tn{sup}}
			\left[
			t(1-t)\|\tilde{\nu}_{k}-\tilde{\nu}_{l}\|_{\tilde{\bSig}_{kl}^{-1}(t)}^{2}
			\right]}.
	\end{equation}
	Here we have suppressed the functional dependence on the underlying model parameters $\bB$ and $\bPi$.
	For large $n$, observe that as $\rho_{\tn{A}}^{\star}$ increases, $\rho_{\tn{A}}$ also increases, and therefore the large-sample optimal error rate corresponding to adjacency spectral embedding decreases in light of Eq.~(\ref{eq:limitErrorChernoff}) and its generalization. Similarly, large values of $\rho_{\tn{L}}^{\star}$ correspond to good theoretical performance of Laplacian spectral embedding. Thus, if $\rho^{\star} > 1$, then ASE is to be preferred to LSE, whereas if $\rho^{\star} < 1$, then LSE is to be preferred to ASE. The case when $\rho^{\star} = 1$ indicates that neither ASE nor LSE is superior for the given parameters $\bB$ and $\bPi$. To reiterate, we summarize these preferences as ASE~$>$~LSE, ASE~$<$~LSE, and ASE~$=$~LSE, respectively.
	
	In what follows, we fixate on the asymptotic quantity $\rho^{\star}$. 
	For the two-block SBM and certain $K$-block SBMs exhibiting symmetry, Eq.~(\ref{eq:rhoStar}) reduces to the simpler form
	\begin{equation}
	\label{eq:rhoStarSimple}
		\rho^{\star}
		= \frac{\underset{t\in(0,1)}{\tn{sup}}
			\left[
			t(1-t)\|\nu_{1}-\nu_{2}\|_{\bSig_{1,2}^{-1}(t)}^{2}\right]}{
		\underset{t\in(0,1)}{\tn{sup}}
			\left[
			t(1-t)\|\tilde{\nu}_{1}-\tilde{\nu}_{2}\|_{\tilde{\bSig}_{1,2}^{-1}(t)}^{2}
			\right]}
	\end{equation}	
	for canonically specified latent positions $\nu_{1}$ and $\nu_{2}$. In some cases it is possible to concisely obtain analytic expressions (in $t$) for both the numerator and denominator. In other cases this is not possible. A related challenge with respect to Eq.~(\ref{eq:rhoStar}) is analytically inverting the interpolated block conditional covariance matrices $\bSig_{1,2}(t)$ and $\tilde{\bSig}_{1,2}(t)$. Section~\ref{sec:Appendix} provides additional technical details and discussion addressing these issues.
	
	%For the \authorER\ regime, the notion of measuring relative performance via Chernoff information is meaningless since block membership is known vacuously for all vertices (i.e.~there is effectively a single block). In this case, $\rho^{\star}$ exhibits a singularity, although it would be logical to define $\rho^{\star} := 1$ in light of relative spectral performance being moot.	
	% XXXXXXXXXXXXXXXXXXXXXXXXXXXXXXXXXXXXXXXXXXXXXXXXXXXXXXXXXXXXXXXXXXXXXXXX
	\section{Elucidating network structure}
	\label{sec:ElucidatingNetworkStructure}
	\subsection{The two-block stochastic block model}
	Consider the set of two-block SBMs with parameters $\bPi\equiv(\pi_{1},1-\pi_{1})$
	and
	$\bB \in \mathscr{B}:=\left\{
		\bB = \Bigl[
			\begin{smallmatrix}
				a & b \\
				b & c
			\end{smallmatrix}
		\Bigr]: a,b,c \in (0,1)\right\}$.
	For $\bPi=(\tfrac{1}{2},\tfrac{1}{2})$, then $a \ge c$ without loss of generality by symmetry.
	In general, for any fixed choice of $\bPi$, the class of models $\mathscr{B}$ can be
	%identified by the open unit cube, $(0,1)^{3}$, up to the symmetry $a \ge c$. Moreover, $\mathscr{B}$ can be
	partitioned according to matrix rank, namely
	\begin{align*}
		\label{eq:scB_rank_decomp}
		\mathscr{B}
		&\equiv \mathscr{B}_{1} \bigsqcup \mathscr{B}_{2}
		:= \{\bB:\tn{rank}(\bB)=1; a,b,c \in (0,1)\}
		\bigsqcup\{\bB:\tn{rank}(\bB)=2; a,b,c \in (0,1)\}.
	\end{align*}
	The collection of sub-models $\mathscr{B}_{1}$ further decomposes into the disjoint union of
	the Erd\H{o}s--R\'{e}nyi model with homogeneous edge probability $a=b=c\in(0,1)$ and its relative complement in $\mathscr{B}_{1}$ satisfying the determinant constraint $\rmDet(\bB) \equiv ac-b^{2} = 0$. These partial sub-models can be viewed as one-dimensional and two-dimensional (parameter) regions in the open unit cube, $(0,1)^{3}$, respectively.
	
	Similarly, the collection of sub-models $\mathscr{B}_{2}$ further decomposes into the disjoint union of $\mathbb{PD}_{2}\cap\mathscr{B}_{2}$ and $\mathbb{IND}_{2}\cap\mathscr{B}_{2}$, where $\mathbb{PD}_{2}$ denotes the set of positive definite matrices in $\R^{2 \times 2}$ and $\mathbb{IND}_{2} :=\{\bB\in\mathscr{B}_{2}: \exists \bX\in\R^{2 \times 2}, \rmRank(\bX)=2, \bB=\bX\bI_{1}^{1}\bX^{\top}\}$.
	%or $\textnormal{rank}(\bB)=2$ and since $a \ge c$, then $\bB \in \mathbb{PD}_{2}$ if and only if both $a>b$ and $ac > b^{2}$. Otherwise, $\bB\notin\mathbb{PD}_{2}$ which holds if and only if either $a < b$ (hence $ac \le a^{2} < b^{2}$) or collectively $a \ge b$ and $ ac < b^{2}$.
	Here only $\bI_{0}^{2}\equiv\bI_{2}$ and $\bI_{1}^{1}$ are necessary for computing edge probabilities via inner products of the latent positions. Both of these partial sub-models can be viewed as three-dimensional (parameter) regions in $(0,1)^{3}$.
	
	\begin{remark}[Latent position parametrization]
	\label{rem:LatentPositionParam}
		One might ask whether or not for our purposes there exists a ``best'' latent position representation for some or even every SBM. To this end and more generally, for any $K \ge 2$ and $\bM\in\mathbb{PD}_{K} \subset \R^{K \times K}$, there exists a unique lower-triangular matrix $\bL\in\R^{K \times K}$ with positive diagonal entries such that $\bM=\bL\bL^{\top}$ by the Cholesky matrix decomposition. This yields a canonical choice for the matrix of latent positions $\bX$ when $\bB$ is positive definite. In particular, for $\bB\in\mathbb{PD}_{2}$, then $\bB=\bX\bI_{2}\bX^{\top}$ with
	$\bX := \Bigl[
		\begin{smallmatrix}
			\sqrt{a} & 0 \\
			b/\sqrt{a} & \sqrt{ac-b^{2}}/\sqrt{a}
		\end{smallmatrix}
	\Bigr]$.
	In contrast, for $\bB\in\mathbb{IND}_{2}$, then $\bB = \bX\bI^{1}_{1}\bX^{\top}$ with
	$\bX := \Bigl[
		\begin{smallmatrix}
			\sqrt{a} & 0 \\
			b/\sqrt{a} & \sqrt{b^{2}-ac}/\sqrt{a}
		\end{smallmatrix}
	\Bigr]$,
	keeping in mind that in this case $b^{2}-ac > 0$. The latter factorization may be viewed informally as an 
	indefinite Cholesky decomposition under $\bI_{1}^{1}$. For the collection of rank one sub-models $\mathscr{B}_{1}$, the latent positions $\nu_{1}$ and $\nu_{2}$ are simply taken to be scalar-valued.
	\end{remark}

	\subsubsection{Homogeneous balanced network structure}
	\label{sec:K2_homo_balanced}

	We refer to the two-block SBM sub-model with
	$\bB = \Bigl[\begin{smallmatrix}
				a & b \\
				b & a
	\end{smallmatrix}\Bigr]$
	and $\bPi = (\tfrac{1}{2},\tfrac{1}{2})$
	as the \emph{homogeneous balanced two-block SBM}. The cases when $a>b$, $a<b$, and $a=b$ correspond to the cases when $\bB$ is positive definite, indefinite, and reduces to \authorER, respectively.
	The positive definite parameter regime has the network structure interpretation of being \emph{assortative} in the sense that the within-block edge probability $a$ is larger than the between-block edge probability $b$, consistent with the affinity-based notion of community structure. In contrast, the indefinite parameter regime has the network structure interpretation of being \emph{disassortative} in the sense that between-block edge density exceeds within-block edge density, consistent with the ``opposites attract'' notion of community structure.
	
	For this SBM sub-model, $\rho^{\star}$ can be simplified analytically (see Section~\ref{sec:Appendix} for additional details) and can be expressed as a translation with respect to the value one, namely
	\begin{equation}
		\label{eq:K2rhoStar}
		\rho^{\star}
		\equiv \rho^{\star}_{a,b}
		= 1 + \frac{(a-b)^{2}(3a(a-1)+3b(b-1)+8ab)}{4(a+b)^{2}(a(1-a)+b(1-b))}
		:= 1 + c_{a,b} \times \psi_{a,b},
	\end{equation}
	where $\psi_{a,b}:=3a(a-1)+3b(b-1)+8ab$ and $c_{a,b}>0$.
	By recognizing that $\psi_{a,b}$ functions as a discriminating term, it is straightforward to read off the relative performance of ASE and LSE according to Table~\ref{table:K2_homo_relative_performance}.

	\begin{figure}
		\centering
		\includegraphics[width=0.5\textwidth]{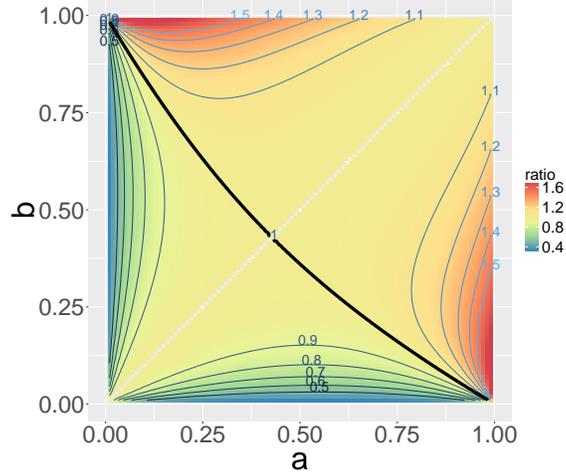}
		\caption{The ratio $\rho^{\star}$ for the homogeneous balanced sub-model in Section \ref{sec:K2_homo_balanced}. The empty diagonal depicts the \authorER\ model singularity at $a=b$.}
		\label{fig:plot_ssbmAB}
	\end{figure}

	\begin{table}[h]
		\caption{Summary of embedding performance in Section~\ref{sec:K2_homo_balanced}}
		\label{table:K2_homo_relative_performance}
	\begin{center}
		\begin{tabular}{c}
			\hline 
			$\rho^{\star}=1 \Longleftrightarrow \psi_{a,b}=0$ \textnormal{; (ASE~$=$~LSE)}\\
			$\rho^{\star}>1 \Longleftrightarrow \psi_{a,b}>0$ \textnormal{; (ASE~$>$~LSE)}\\
			$\rho^{\star}<1 \Longleftrightarrow \psi_{a,b}<0$ \textnormal{; (ASE~$<$~LSE)}
		\end{tabular}
	\end{center}
	\end{table}
	
	Further investigation of Eq.~(\ref{eq:K2rhoStar}) leads to the observation that ASE~$<$~LSE for all $0<b<a \le \tfrac{3}{7}$, thereby yielding a parameter region for which LSE dominates ASE. On the other hand, for any fixed $b\in(0,1)$ there exist values $a_{1} < a_{2}$ such that ASE~$<$~LSE under $a_{1}$, whereas ASE~$>$~LSE under $a_{2}$. Figure~\ref{fig:plot_ssbmAB} demonstrates that for homogeneous balanced network structure, LSE is preferred to ASE when the entries in $\bB$ are sufficiently small, whereas conversely ASE is preferred to LSE when the entries in $\bB$ are not too small.
	\begin{remark}[Model spectrum and ASE dominance I]
	\label{rem:ModelSpectrumASE_I}
		In the current setting $\lambda_{\textrm{max}}(\bB)=a+b$, hence $\lambda_{\textrm{max}}(\bB)>1$ implies ASE~$>$~LSE by Eq.~(\ref{eq:K2rhoStar}). This observation amounts to a network structure-based (i.e.~$\bB$-based) spectral sufficient condition for determining when ASE is preferred to LSE.
	\end{remark}
	\begin{remark}[A balanced one-dimensional SBM restricted sub-model]
	When $b=1-a$, the homogeneous balanced sub-model further reduces to a one-dimensional parameter space such that $\rho^{\star}$ simplifies to
	\begin{equation} 
		\rho^{\star} = 1+\tfrac{1}{4}(2a-1)^2 \ge 1,
	\end{equation}
	demonstrating that ASE uniformly dominates LSE for this restricted sub-model. Additionally, it is potentially of interest to note that in this setting the marginal covariance matrices from Theorem~\ref{thrm:GRDPG_CLT_ASE} for ASE coincide for each block. In contrast, the same behavior is not true for LSE.
	\end{remark}
	
	\subsubsection{Core-periphery network structure}
	\label{sec:K2_core_periphery}
	
%	\begin{figure}[h]
%		\centering
%		\includegraphics[width=0.4\textwidth]{PLOTcp.pdf}
%		\caption{The ratio $\rho^{\star}$ for the core-periphery balanced sub-model in Section \ref{sec:K2_core_periphery}. The empty diagonal depicts the \authorER\ model singularity at $a=b$.}
%		\label{fig:plot_cp}
%	\end{figure}
	
	%previously \centering\includegraphics[width=0.8\textwidth]{PLOTcp050.pdf}

	We refer to the two-block SBM sub-model with
	$\bB = \Bigl[\begin{smallmatrix}
		a & b \\
		b & b
	\end{smallmatrix}\Bigr]$
	and $\bPi = (\pi_{1}, 1-\pi_{1})$
	as the \emph{core-periphery two-block SBM}. We explicitly consider the balanced (block size) regime in which $\bPi = (\tfrac{1}{2}, \tfrac{1}{2})$ and an unbalanced regime in which $\bPi = (\tfrac{1}{4}, \tfrac{3}{4})$. Here, the cases $a>b$, $a<b$, and $a=b$ correspond to the cases when $\bB$ is positive definite, indefinite, and reduces to the \authorER\ model, respectively.
	
\begin{figure}
		\centering
		\begin{subfigure}[b]{0.5\linewidth}
			\centering\includegraphics[width=\textwidth]{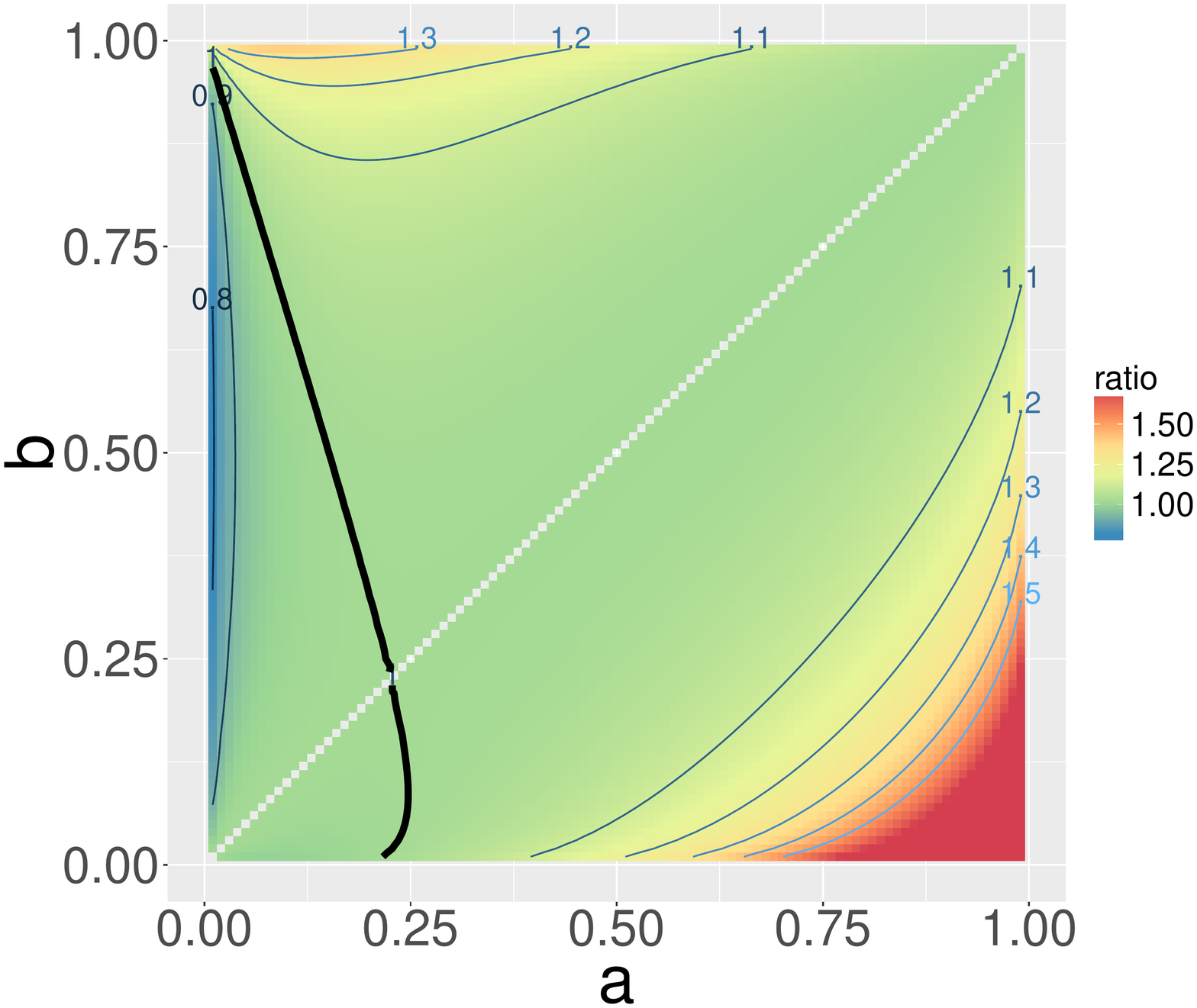}
			\caption{\label{fig:plot_cp050} $\rho^{\star}$ when $\bPi=(\tfrac{1}{2},\tfrac{1}{2})$ i.e.~balanced}
		\end{subfigure}%
		\begin{subfigure}[b]{0.5\linewidth}
			\centering\includegraphics[width=\textwidth]{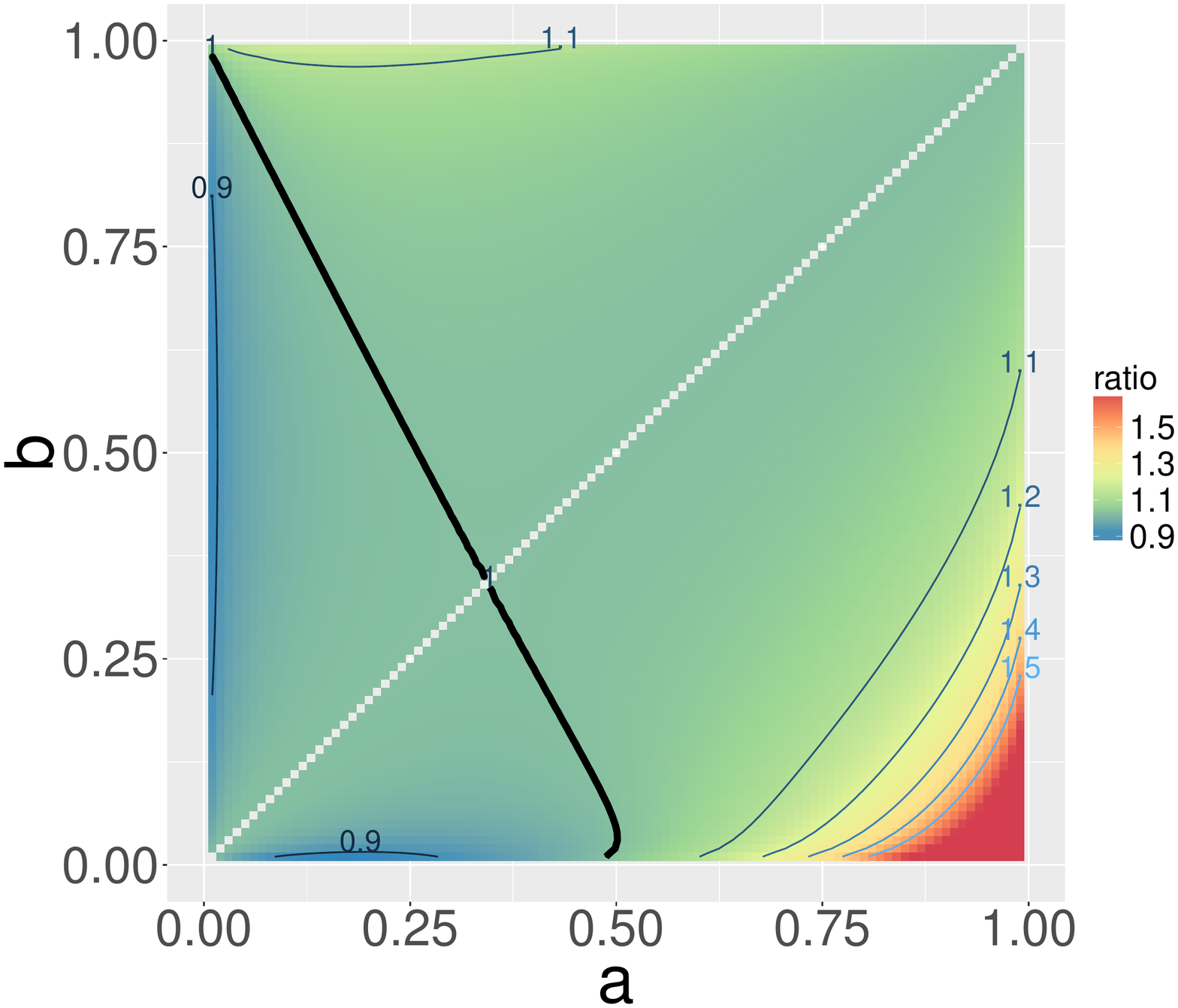}
			\caption{\label{fig:plot_cp025} $\rho^{\star}$ when $\bPi=(\tfrac{1}{4},\tfrac{3}{4})$ i.e.~unbalanced}
		\end{subfigure}
		\caption{\label{fig:plot_cp}The ratio $\rho^{\star}$ for the core-periphery sub-model in Section \ref{sec:K2_core_periphery}. The empty diagonal depicts the \authorER\ model singularity at $a=b$.}
	\end{figure}

	For this sub-model, the ratio $\rho^{\star}$ is not analytically tractable in general. That is to say, simple closed-form solutions do not simultaneously exist for the numerator and denominator in the definition of $\rho^{\star}$. As such, Figure~\ref{fig:plot_cp} is obtained numerically by evaluating $\rho^{\star}$ on a grid of points in $(0,1)^{2}$ followed by smoothing.
	
	For $a>b$, graphs generated from this SBM sub-model exhibit the popular interpretation of core-periphery structure in which vertices forming a dense core are attached to surrounding periphery vertices with comparatively smaller edge connectivity. Provided the core is sufficiently dense, namely for $a>\tfrac{1}{4}$ in the balanced regime and $a > \tfrac{1}{2}$ in the unbalanced regime, Figure~\ref{fig:plot_cp} demonstrates that ASE~$>$~LSE. Conversely, ASE~$<$~LSE uniformly in $0<b<a$ for small enough values of $a$ in both the balanced and unbalanced regime.
	
	In contrast, when $a<b$, the sub-model produces graphs whose network structure is interpreted as having a comparatively sparse induced subgraph which is strongly connected to all vertices in the graph but for which the subgraph vertices exhibit comparatively weaker connectivity. Alternatively, the second block may itself be viewed as a dense core which is simultaneously densely connected to all vertices in the graph. Figure~\ref{fig:plot_cp} illustrates that for the balanced regime, LSE is preferred for sparser induced subgraphs. Put differently, for large enough dense core with dense periphery, then ASE is the preferable spectral embedding procedure.
	LSE is preferred to ASE in only a  relatively small region corresponding approximately to the triangular region where $0< b < 1 - 4a$, which as a subset of the unit square has area $\tfrac{1}{8}$. Similar behavior holds for the unbalanced regime for approximately the (enlarged) triangular region of the parameter space where $0 < b < 1-2a$, which as a subset of the unit square has area $\tfrac{1}{4}$.
	
	Figure~\ref{fig:plot_cp} suggests that as $\pi_{1}$ decreases from $\tfrac{1}{2}$ to $\tfrac{1}{4}$, LSE is favored in a growing region of the parameter space, albeit still in a smaller region than that for which ASE is to be preferred. Together with the observation that LSE dominates in the lower-left corner of the plots in Figure~\ref{fig:plot_cp} where $a$ and $b$ have small magnitude, we are led to say in summary that LSE favors relatively sparse core-periphery network structure. To reiterate, sparsity is interpreted with respect to the parameters $a$ and $b$, keeping in mind the underlying simplifying assumption that $n_{k} = n \pi_{k}$ for $k = 1, 2$.

	\begin{remark}[Model spectrum and ASE dominance II]
	\label{rem:ModelSpectrumASE_II}	
		For $0<b<a<1$, then $\lambda_{\textrm{max}}(\bB)=\tfrac{1}{2}\left(a+b+\sqrt{a^{2}-2ab+5b^{2}}\right)$. Numerical evaluation (not shown) yields that $\lambda_{\textrm{max}}(\bB)>\tfrac{1}{2}$ implies ASE~$>$~LSE. Along the same lines as the discussion in Section~\ref{sec:K2_homo_balanced}, this observation provides a network structure (i.e.~$\bB$-based) spectral sufficient condition for this sub-model for determining the relative embedding performance ASE~$>$~LSE.
	\end{remark}
	
	\subsubsection{Two-block rank one sub-model}
	\label{sec:K2_rank1_submodel}
	
	The sub-model for which
	$\bB = \Bigl[\begin{smallmatrix}
				a & b \\
				b & c
			\end{smallmatrix}\Bigr]$
	with $a, b, c \in (0,1)$ and $\rmDet(\bB)=0$ can be re-parameterized according to the assignments $a \mapsto p^{2}$ and $c \mapsto q^{2}$, yielding
	$\bB = \Bigl[\begin{smallmatrix}
				p^{2} & pq \\
				pq & q^{2}
			\end{smallmatrix}\Bigr]$
	with $p, q\in(0,1)$. Here $\tn{rank}(\bB)=1$ and $\bB$ is positive semidefinite, corresponding to the one-dimensional RDPG model with latent positions given by the scalars $p$ and $q$ with associated probabilities $\pi_{1}$ and $\pi_{2}$, respectively. Explicit computation yields the expression
		\begin{equation}
			\label{eq:ratio_rank1_analytic}
			\fontsize{10}{0}
			\rho^{\star} =
			\frac{(\sqrt{p}+\sqrt{q})^{2}(\pi_{1}p^{2}+\pi_{2}q^{2})^{2}\left(\sqrt{\pi_{1} p(1-p^{2})+\pi_{2}q(1-pq)}+\sqrt{\pi_{1} p(1-pq)+\pi_{2}q(1-q^{2})}\right)^{2}}{4(\pi_{1}p+\pi_{2}q)^{2}\left(\sqrt{\pi_{1}p^{4}(1-p^{2})+\pi_{2}p q^{3}(1-pq)}+\sqrt{\pi_{1}p^{3}q(1-pq)+\pi_{2}q^{4}(1-q^{2})}\right)^{2}},
		\end{equation}
	whereby $\rho^{\star}$ is given as an explicit, closed-form function of the parameter values $p$, $q$, and $\pi_{1}$ with $\pi_{2}=1-\pi_{1}$. The simplicity of this sub-model together with its analytic tractability with respect to both $\bB$ and $\bPi$ makes it particularly amenable to study for the purpose of elucidating network structure. Below, consideration of this sub-model further illustrates the relationship between (parameter-based) sparsity and relative embedding performance.
	%To reiterate, undergirding our analysis is the simplifying assumption that $n_{1}=\pi_{1}n$ and $n_{2}=\pi_{2}n$ for $n$-vertex SBM graphs.
	
	\begin{figure}
		\centering
		\begin{subfigure}[b]{0.5\linewidth}
			\centering\includegraphics[width=\textwidth]{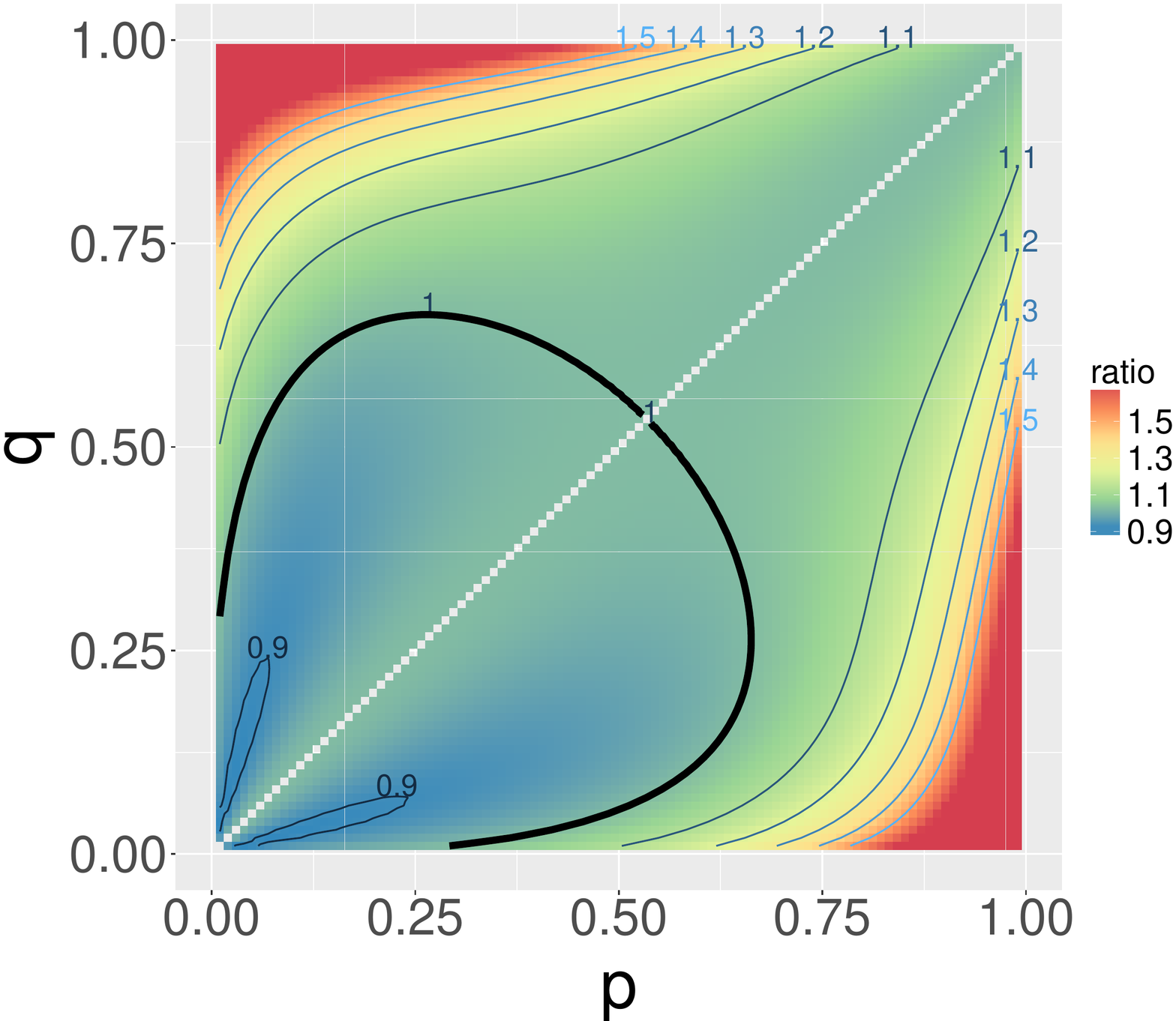}
			\caption{\label{fig:plot_rank1pq_pi1/2} $\rho^{\star}$ when $\bPi=(\tfrac{1}{2},\tfrac{1}{2})$ i.e.~balanced}
		\end{subfigure}%
		\begin{subfigure}[b]{0.5\linewidth}
			\centering\includegraphics[width=\textwidth]{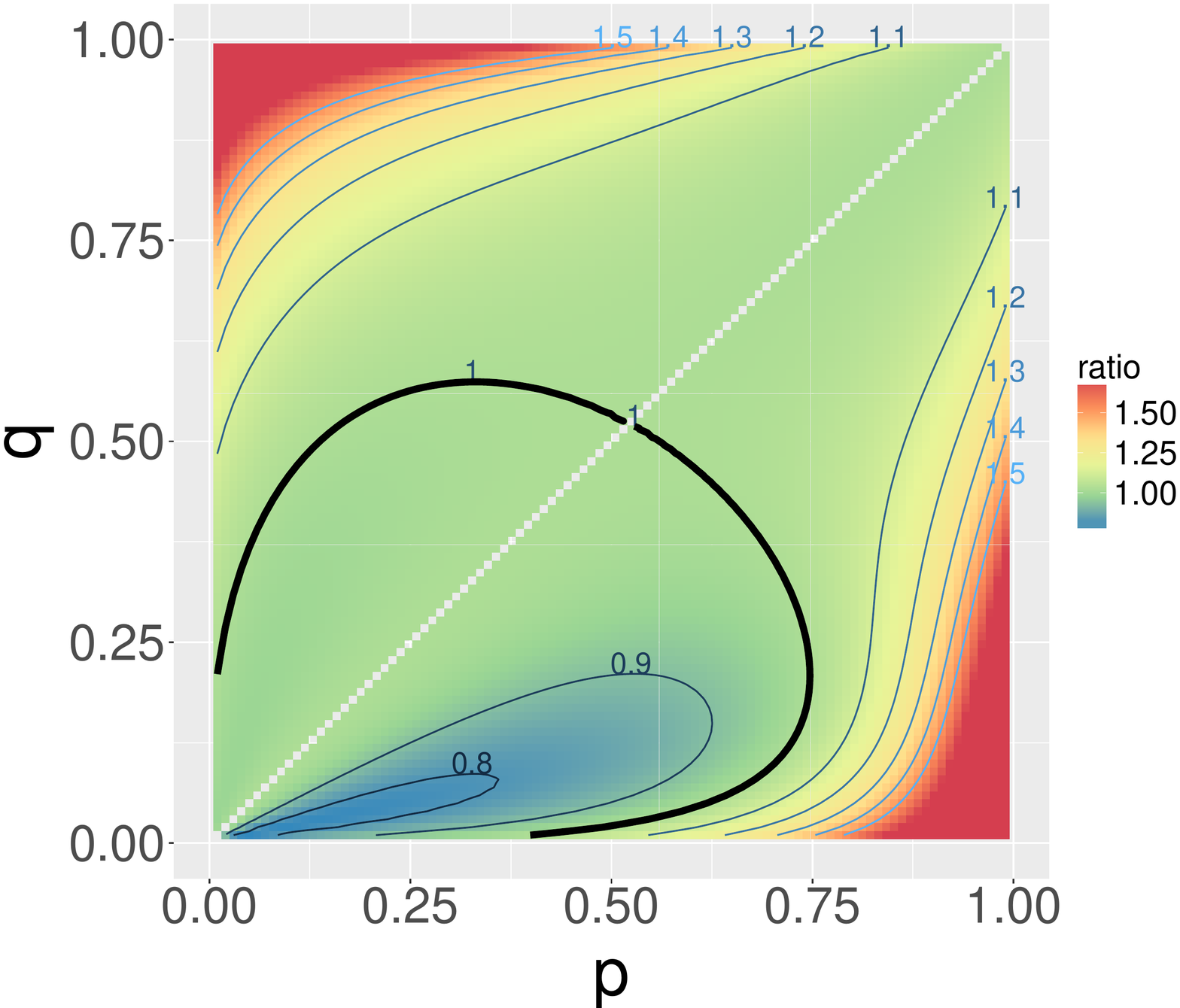}
			\caption{\label{fig:plot_rank1pq_pi1/4} $\rho^{\star}$ when $\bPi=(\tfrac{1}{4},\tfrac{3}{4})$ i.e.~unbalanced}
		\end{subfigure}
		\caption{\label{fig:plot_rank1pq}The ratio $\rho^{\star}$ for the two-block rank one sub-model in Section~\ref{sec:K2_rank1_submodel}. The empty diagonal depicts the \authorER\ model singularity at $p=q$.}
	\end{figure}
	
	Figure~\ref{fig:plot_rank1pq} demonstrates how LSE favors sparse graphs in the sense of the edge probabilities, $p$ and $q$, as well as how relative performance changes in light of (un)balanced block sizes, reflected by $\pi_{1}$. Here the underlying $\bB$ matrix is always positive semidefinite, and each of the regions $p > q$ and $p < q$ corresponds to a modified notion of core-periphery structure. For example, when $p>q$, then
	$\bB =
		\Bigl[\begin{smallmatrix}
			p_{1} & p_{2} \\
			p_{2} & p_{3}
		\end{smallmatrix}\Bigr]$
	with $p_{1} > p_{2} > p_{3}$, yielding a hierarchy of core-periphery structure when passing from vertices that are both in block one to vertices that are in different blocks and finally to vertices that are both in block two. Note the similar behavior in the bottom-right triangular regions in Figure~\ref{fig:plot_rank1pq_pi1/2}--\ref{fig:plot_rank1pq_pi1/4} and in the same bottom-right triangular region in Figure~\ref{fig:plot_cp}.

	\begin{remark}[The two-block polynomial $p$ SBM restricted sub-model]
		\begin{figure}
			\centering
			\begin{subfigure}[b]{0.33\linewidth}
				\centering\includegraphics[width=\textwidth]{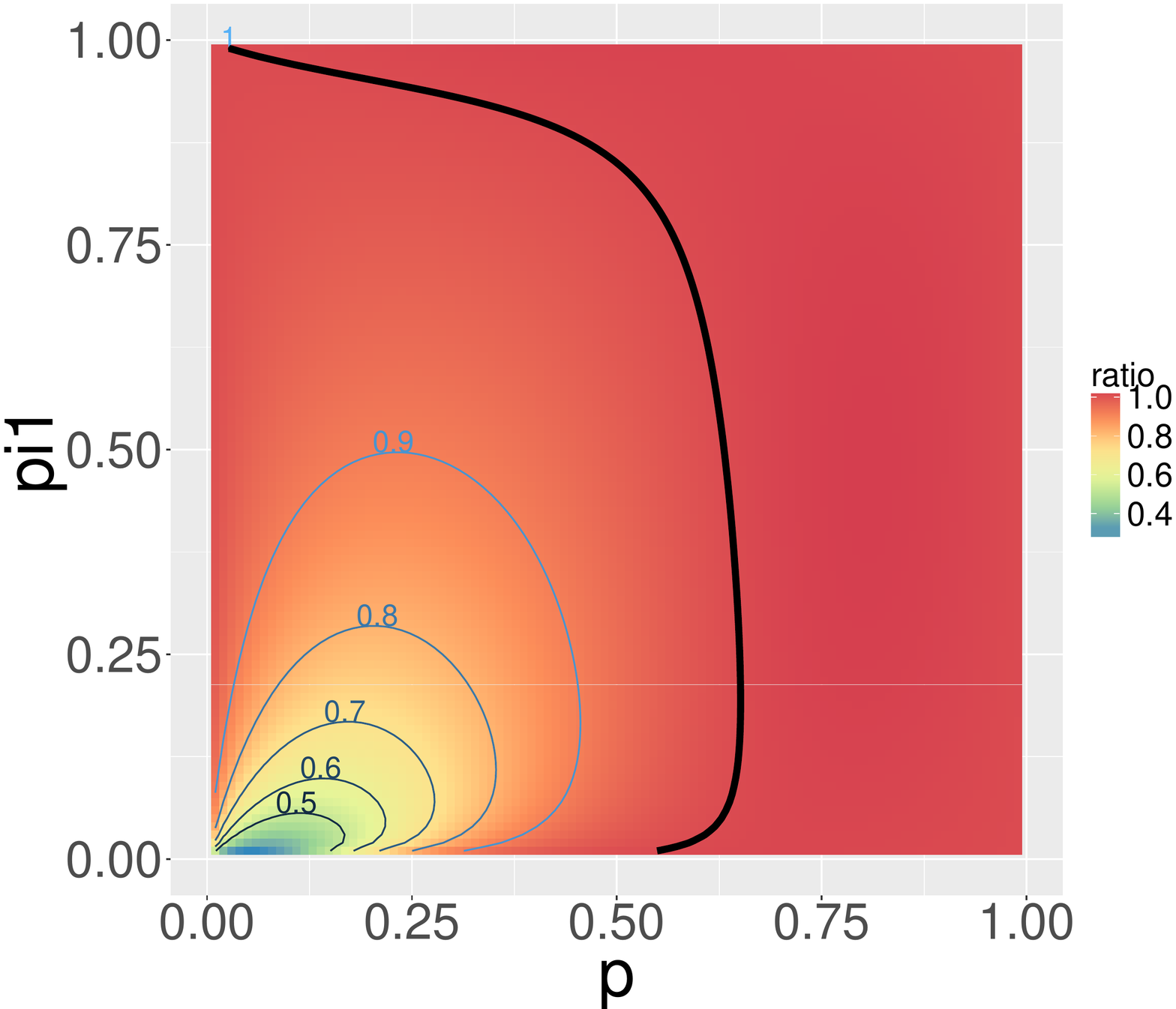}
				\caption{\label{fig:plot_rank1_pp2andpi}
					$\rho^{\star}$ with quadratic $q$
				}
			\end{subfigure}%
			\begin{subfigure}[b]{0.33\linewidth}
				\centering\includegraphics[width=\textwidth]{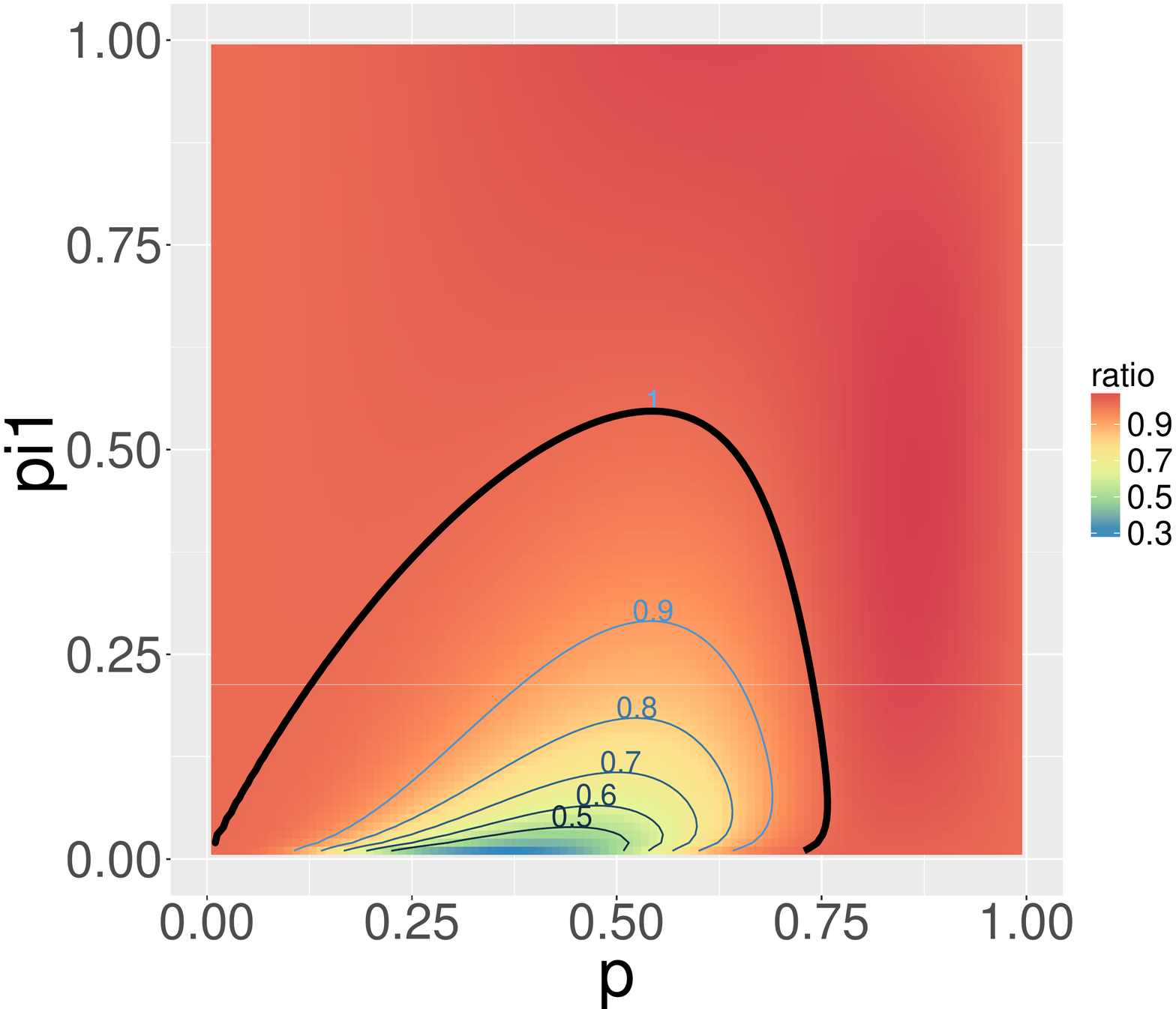}
				\caption{\label{fig:plot_rank1_pp4andpi}
					$\rho^{\star}$ with quartic $q$
				}
			\end{subfigure}
			\begin{subfigure}[b]{0.33\linewidth}
				\centering\includegraphics[width=\textwidth]{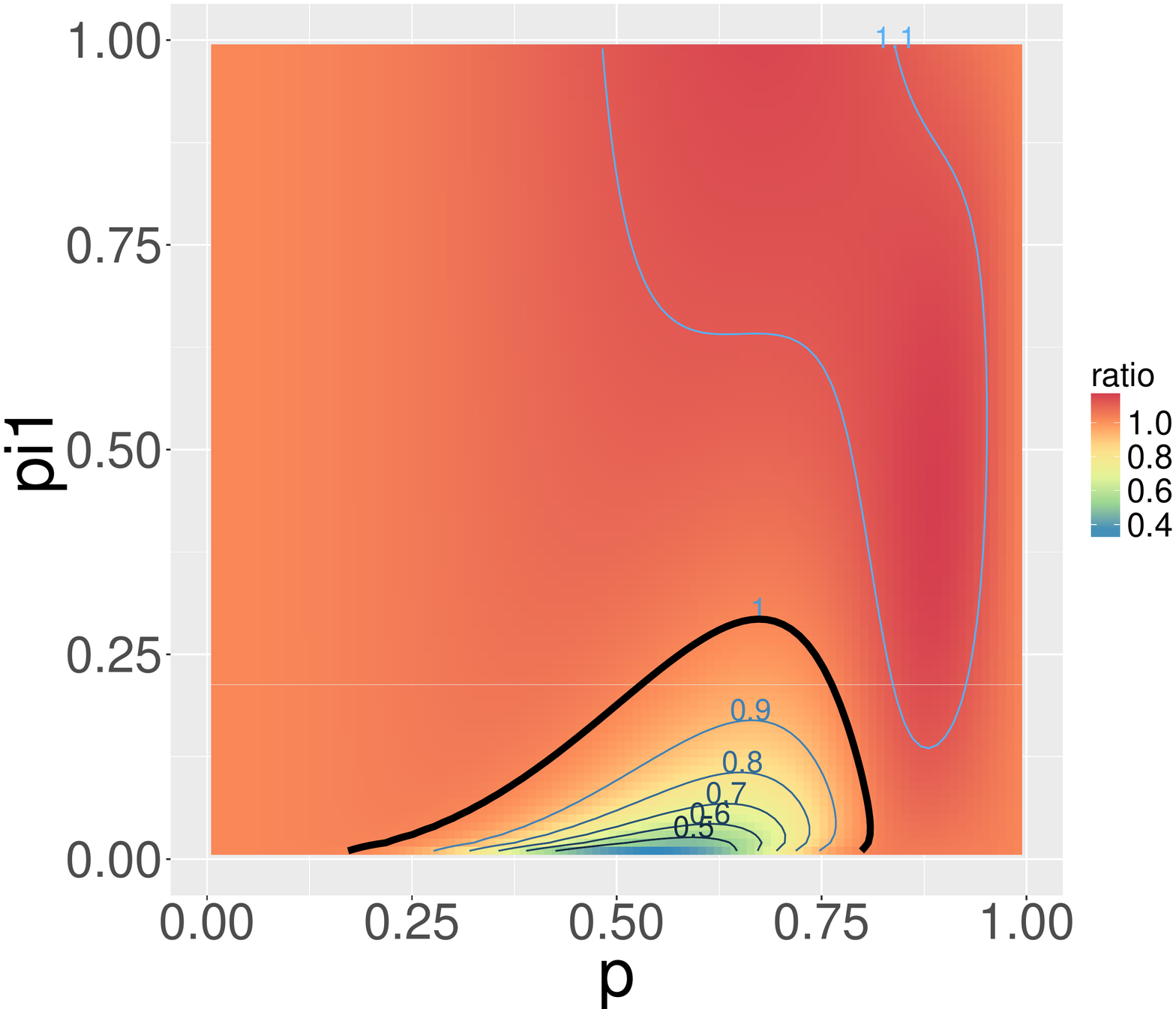}
				\caption{\label{fig:plot_rank1_pp6andpi}
					$\rho^{\star}$ with sextic $q$
				}
			\end{subfigure}%
			\caption{\label{fig:plot_rank1_polypandpi}The ratio $\rho^{\star}$ for $p, \pi_{1} \in (0,1)$ when $q=p^{\gamma}$, $\gamma\in\{2,4,6\}$ in Section~\ref{sec:K2_rank1_submodel}.}
		\end{figure}
		\begin{figure}
			\centering
			\begin{subfigure}[b]{0.45\linewidth}
				\centering\includegraphics[width=\textwidth]{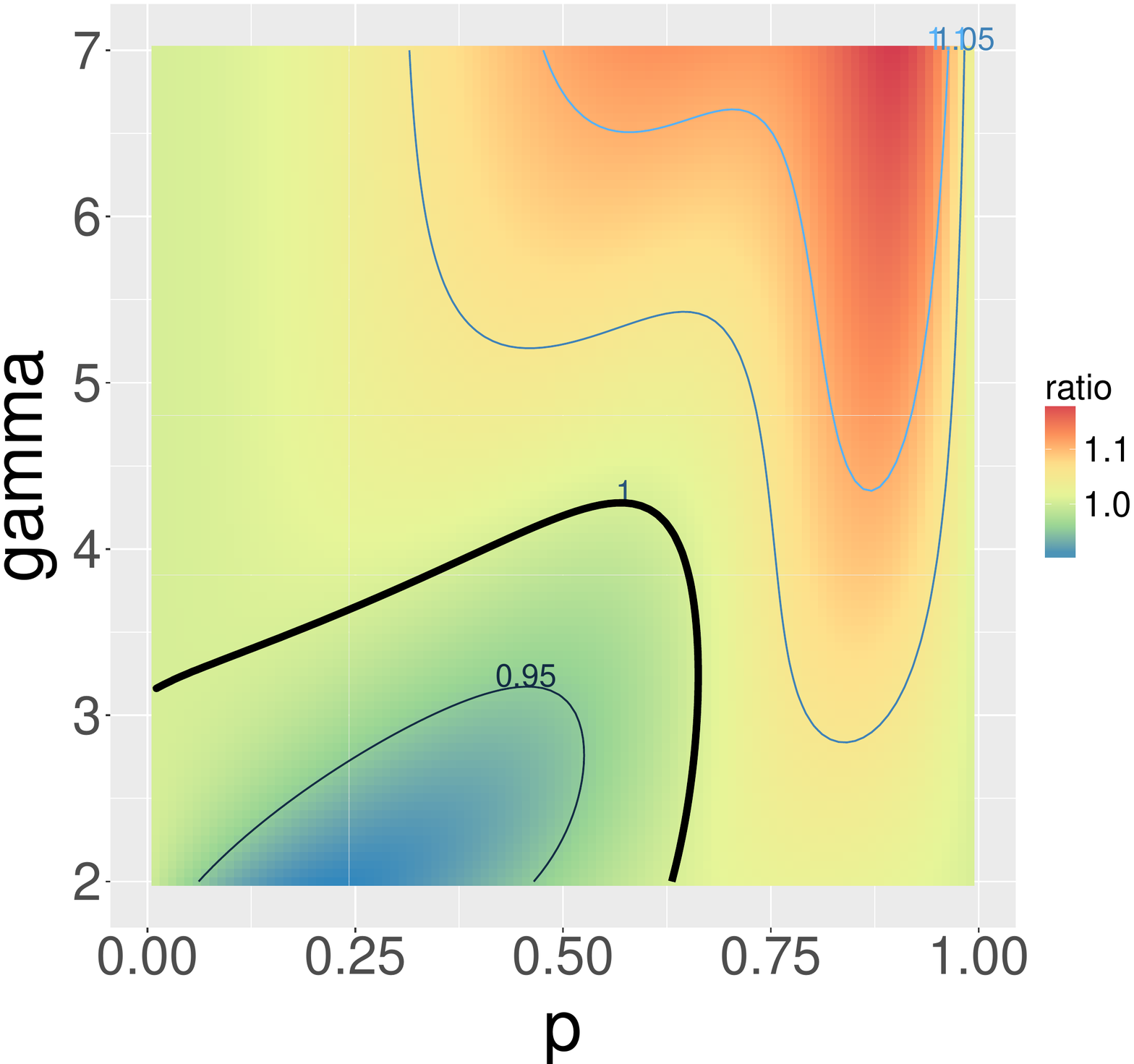}
				\caption{\label{fig:PLOTpolyPgammaPiHalf}
					$\rho^{\star}$ when $\bPi=(\tfrac{1}{2},\tfrac{1}{2})$ i.e.\ balanced
				}
			\end{subfigure}
			\begin{subfigure}[b]{0.45\linewidth}
				\centering\includegraphics[width=\textwidth]{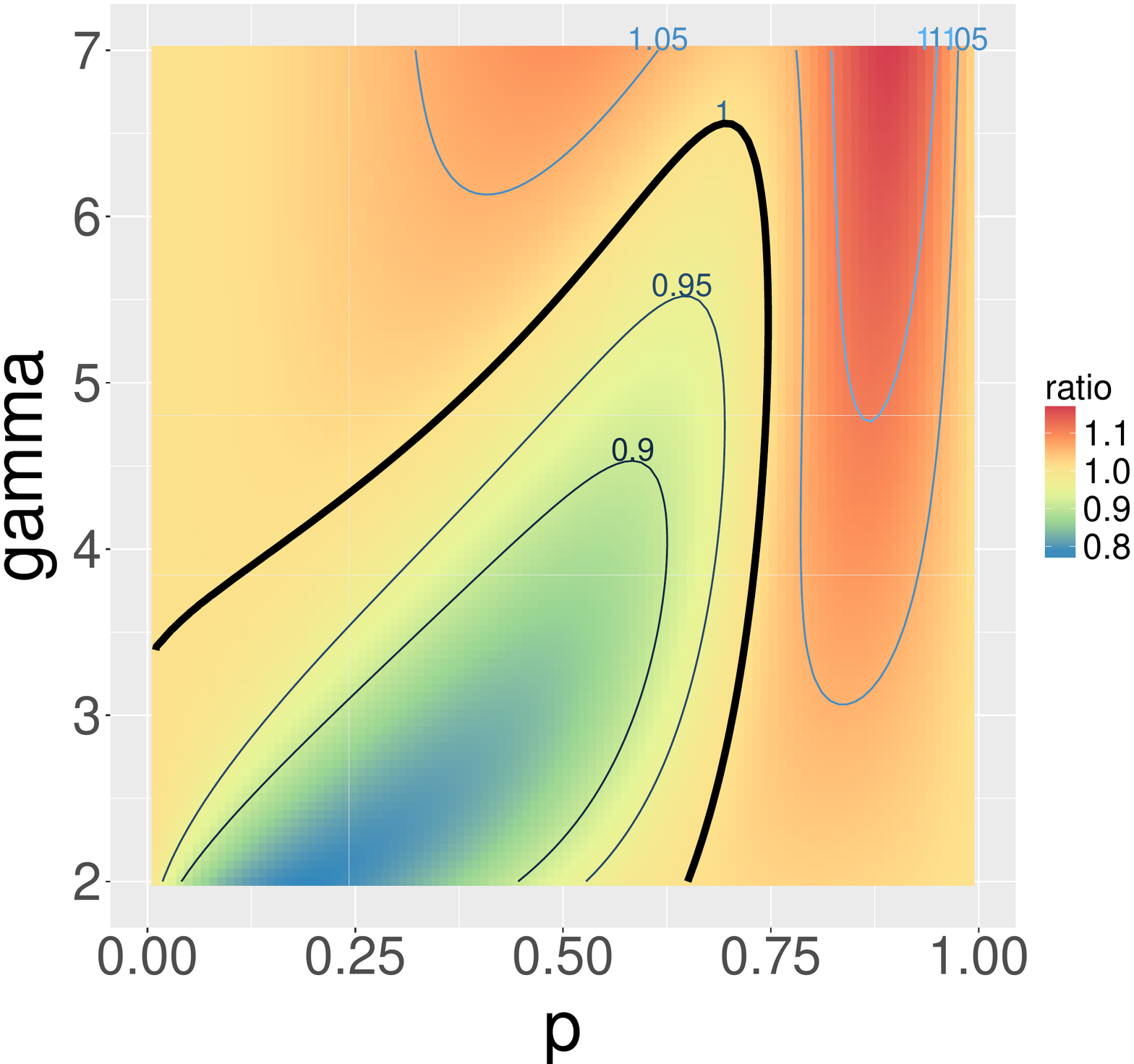}
				\caption{\label{fig:PLOTpolyPgammaPiFourth}
					$\rho^{\star}$ when $\bPi=(\tfrac{1}{4},\tfrac{3}{4})$ i.e.\ unbalanced
				}
			\end{subfigure}%
			\caption{\label{fig:PLOTpolyPgamma}The ratio $\rho^{\star}$ for $p \in (0,1), \gamma \in [2, 7]$ when $q=p^{\gamma}$ in Section~\ref{sec:K2_rank1_submodel}.}
		\end{figure}
		Consider the restricted sub-model in which
		$\bB = \Bigl[\begin{smallmatrix}
					p^{2} & p^{\gamma+1} \\
					p^{\gamma+1} & p^{2\gamma}
				\end{smallmatrix}\Bigr]$,
		where $\gamma > 1$
		and $\pi_{1} \in (0,1)$. For $\gamma \gg 1$ and $\pi_{1}$ fixed, then $\rho^{\star}$ in Eq.~(\ref{eq:ratio_rank1_analytic}) satisfies the approximate behavior
		\begin{equation}
		\label{eq:rhoStar_polyPapprox}
			\rho^{\star}
			\approx \tfrac{\left(1+\sqrt{1-p^{2}}\right)^{2}}{4(1-p^{2})}.
		\end{equation}
		The above approximation exceeds the value one since $1 > \sqrt{1-p^{2}}$ for $p\in(0,1)$ and is simultaneously agnostic with respect to $\pi_{1}$. Moreover, for large values of $\gamma$, the block edge probability matrix is approximately of the form
		$\bB \approx \Bigl[\begin{smallmatrix}
			p_{1} & p_{2} \\
			p_{2} & p_{3}
		\end{smallmatrix}\Bigr]$
		with $p_{1} \gg p_{2} \approx p_{3}$, where $p_{2}$ and $p_{3}$ are very small. This restricted sub-model can therefore be viewed as exhibiting an extremal version of core-periphery structure corresponding to the extremal regions in Figure~\ref{fig:plot_cp} where ASE is preferred.
		
		In Figure~\ref{fig:plot_rank1_polypandpi}, the progression from left to right corresponds to tending towards the approximation presented in Eq.~(\ref{eq:rhoStar_polyPapprox}). For larger values of $\gamma$ when $q=p^{\gamma}$ (not shown), the region where ASE~$>$~LSE continues to expand. We do not discuss or pursue the taking of limits within the parameter space(s) in light of degenerate boundary value behavior and in order to avoid possible misinterpretation.
		
		Figure~\ref{fig:PLOTpolyPgamma} offers a different perspective in which $\gamma$ is allowed to vary continuously for both the balanced and the unbalanced regime. As in Figure~\ref{fig:plot_rank1pq}, Figure~\ref{fig:PLOTpolyPgamma} demonstrates that LSE is preferred for network structure wherein the block with comparatively higher edge probability exhibits smaller block membership size.
	\end{remark}

	\subsubsection{Full rank two-block stochastic block models}
	\label{sec:K2_general_SBM}
		This section presents a macroscopic view of full rank two-block SBMs with
		$\bB =
			\Bigl[\begin{smallmatrix}
				a & b \\
				b & c
			\end{smallmatrix}\Bigr]$,
		$(a,b,c)\in(0,1)^{3}$, for the regimes $\bPi=(\tfrac{1}{2}, \tfrac{1}{2})$ and $\bPi = (\tfrac{1}{4},\tfrac{3}{4})$. The parameter space is partitioned via the latent space geometry of $\bB$, namely according to whether $\bB$ is either positive definite or indefinite.

	\begin{figure}% 3D K=2 plots
		\centering
		\begin{subfigure}[b]{0.5\linewidth}
			\centering\includegraphics[width=\textwidth]{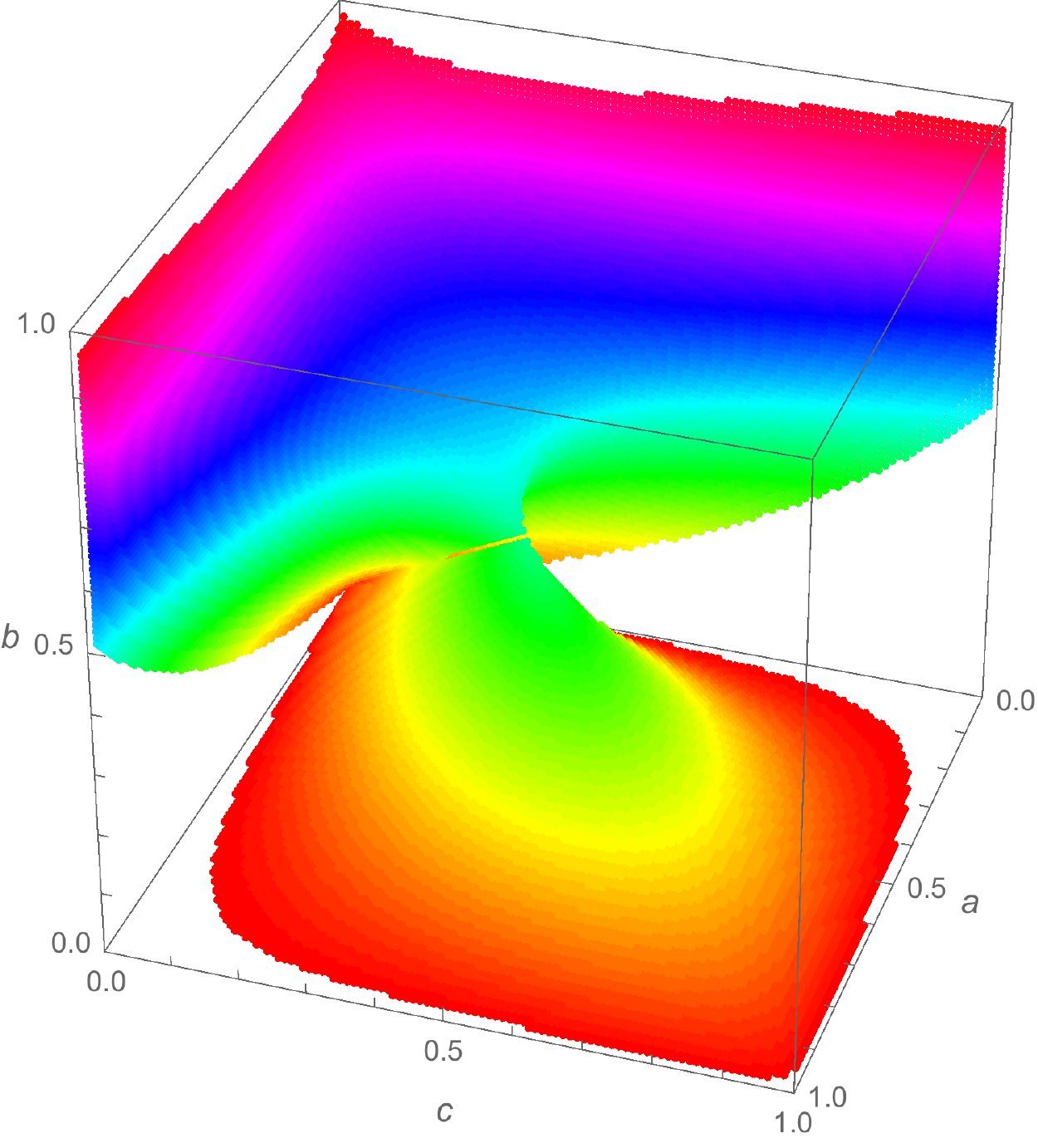}
			\caption{\label{fig:3Dk2pi050} $\rho^{\star}<1$ for $\rmRank(\bB)=2$ when $\bPi=(\tfrac{1}{2},\tfrac{1}{2})$}
		\end{subfigure}%
		\begin{subfigure}[b]{0.5\linewidth}
			\centering\includegraphics[width=\textwidth]{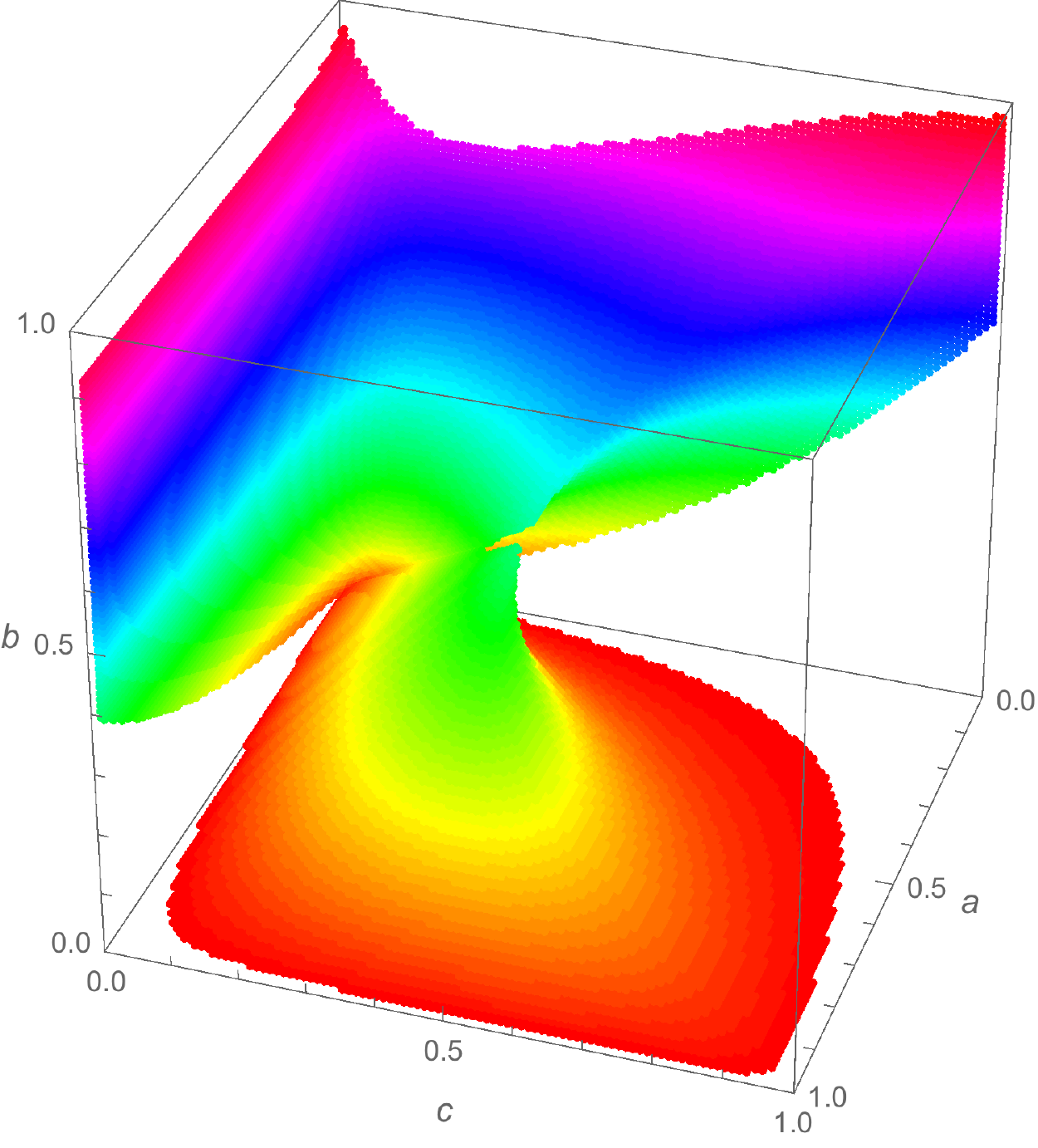}
			\caption{\label{fig:3Dk2pi025} $\rho^{\star}<1$ for $\rmRank(\bB)=2$ when $\bPi=(\tfrac{1}{4},\tfrac{3}{4})$}
		\end{subfigure}
		\caption{\label{fig:3Dk2} The parameter region where ASE~$<$~LSE for full rank $\bB$ in Section~\ref{sec:K2_general_SBM}. The plots depict numerical evaluations of $\rho^{\star}$ for $a,b,c \in [0.01, 0.99]$ with step size 0.01.}
	\end{figure}

	\begin{figure}% 3D top-down PD region plots
		\centering
		\begin{subfigure}[b]{0.5\linewidth}
			\centering
			\includegraphics[width=0.8\textwidth]{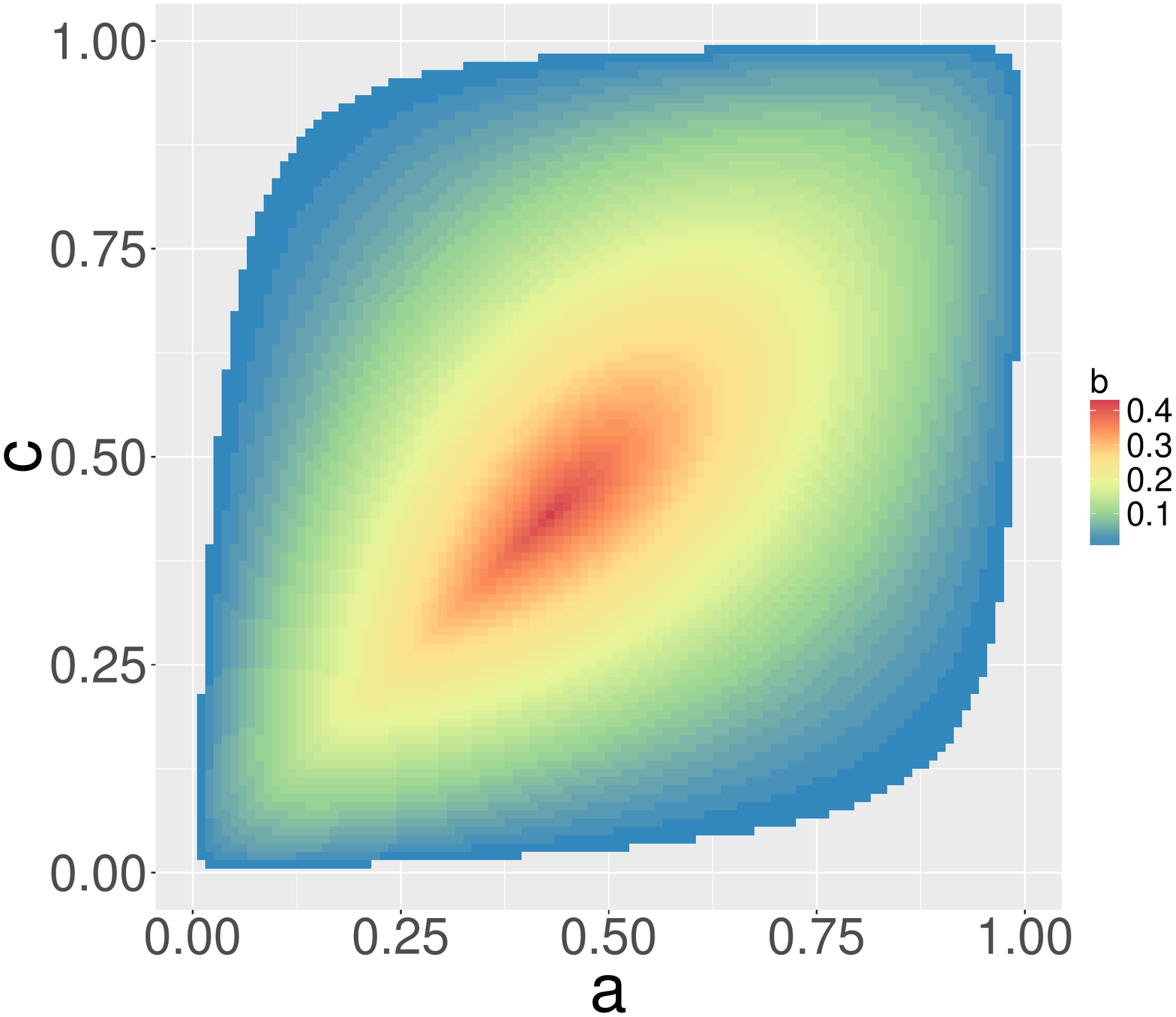}
			\caption{\label{fig:plot_rank2_posdef_pi1/2} $\rho^{\star}<1$ for $\bB \in \mathbb{PD}_{2}$ when $\bPi=(\tfrac{1}{2},\tfrac{1}{2})$}
		\end{subfigure}%
		\begin{subfigure}[b]{0.5\linewidth}
			\centering
			\includegraphics[width=0.8\textwidth]{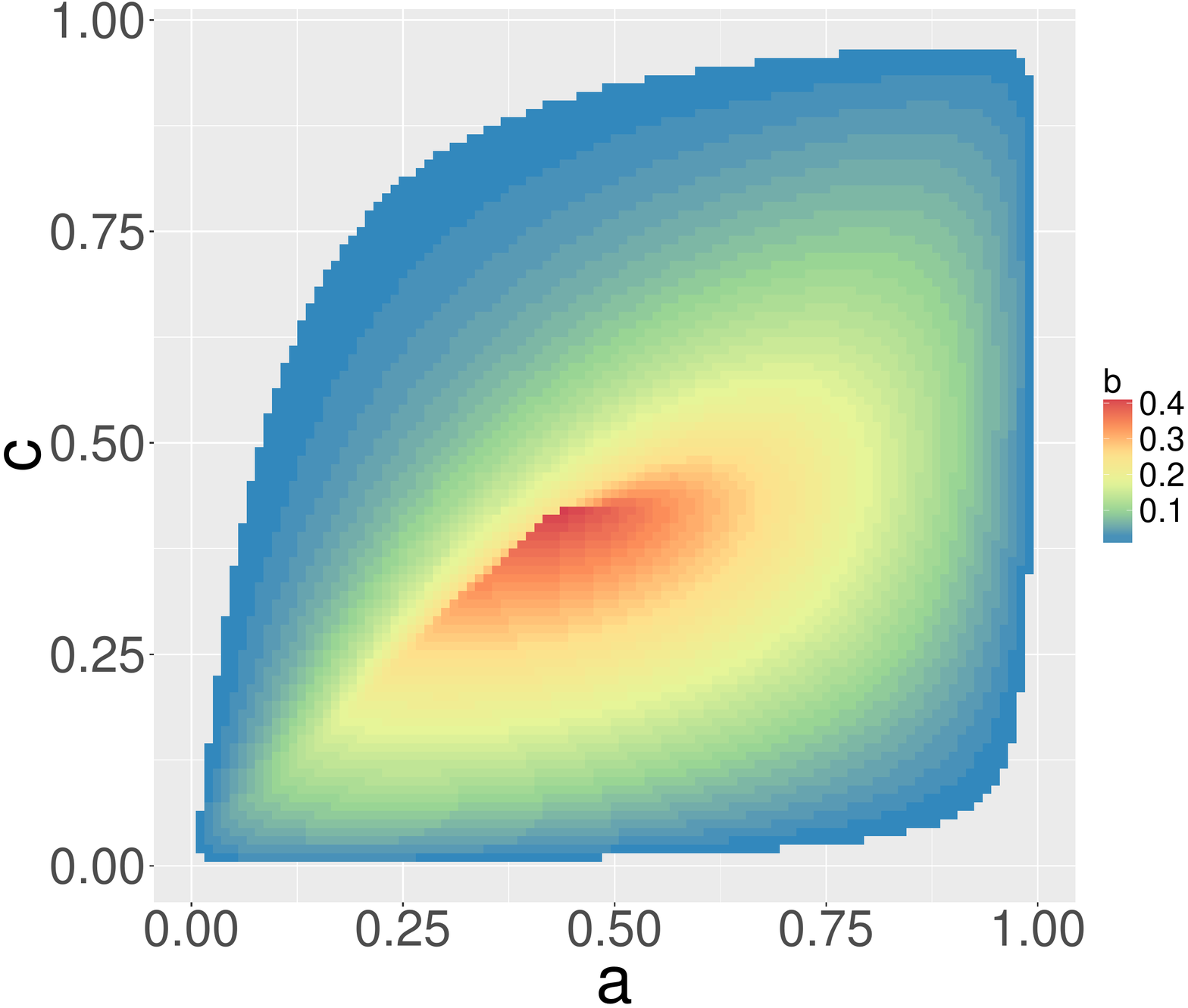}
			\caption{\label{fig:plot_rank2_posdef_pi1/4}$\rho^{\star}<1$ for $\bB \in \mathbb{PD}_{2}$ when  $\bPi=(\tfrac{1}{4},\tfrac{3}{4})$}
		\end{subfigure}
		\caption{\label{fig:plot_rank2_posdef}
			A top-down projected view of the positive definite region where ASE~$<$~LSE in Section~\ref{sec:K2_general_SBM}, with $a$, $b$, and $c$ corresponding to length, depth, and width, respectively. The plots depict numerical evaluations of $\rho^{\star}$ for $a,b,c \in [0.01, 0.99]$ with step size 0.01.}
	\end{figure}

	Figure~\ref{fig:3Dk2pi050} and Figure~\ref{fig:3Dk2pi025} each present a three-dimensional view of the region in the parameter space where ASE~$<$~LSE. The separate positive definite and indefinite parameter regions exhibiting ASE~$<$~LSE can be seen extending from faces of the unit cube. Specifically, the conic-like region rising up from the $b=0$ face corresponds to $\bB$ for which $\bB \in \mathbb{PD}_{2}$, whereas the hyperbolic-like regions extending from the $a=0$ and $c=0$ faces corresponds to $\bB$ for which $\bB \in \mathbb{IND}_{2}$.
	
	For the balanced case reflected in Figure~\ref{fig:3Dk2pi050}, let $a \ge c$ without loss of generality by symmetry, and hence $\rho^{\star}$ is symmetric about the plane defined by $a=c$. For the unbalanced case shown in Figure~\ref{fig:3Dk2pi025}, symmetry no longer holds, and geometric warping behavior can be seen with respect to the $a=c$ plane.
	Figure~\ref{fig:plot_rank2_posdef_pi1/2} and Figure~\ref{fig:plot_rank2_posdef_pi1/4} provide a birds-eye view of the three-dimensional positive definite parameter region from the vantage point $b=``\infty"$. The latter provides another view of the warping phenomenon observed for $\bPi=(\tfrac{1}{4},\tfrac{3}{4})$ that holds in general for all unbalanced regimes.
	
	In both block size regimes depicted in Figure~\ref{fig:3Dk2}, the colored parameter region occupies less than one-fourth of the unit cube volumetrically, thereby quantitatively providing a coarse overall sense in which ASE is to be preferred to LSE for numerous two-block SBM models.

	\subsection{The $K$-block model with homogeneous balanced affinity network structure}
	\label{sec:K_homo_balanced_affinity}
	This section generalizes the analysis in Section~\ref{sec:K2_homo_balanced} to the setting of $K$-block homogeneous balanced affinity SBMs where $\bB_{ij}=a$ for all $i=j$, $\bB_{ij}=b$ for all $i \neq j$, $0<b<a<1$, and $\pi_{i}=\tfrac{1}{K}$ for $1 \le i \le K$.
	\begin{theorem}
		\label{thrm:Kblock_rhoStar}
		For $K$-block homogeneous balanced affinity SBM models as in  Section~\ref{sec:K_homo_balanced_affinity}, the ratio $\rho^{\star}$ in Eq.~(\ref{eq:rhoStar}) can be expressed analytically as
		\begin{equation}
			\label{eq:rhoStar_K_homo_balanced_affinity}
			\rho^{\star}
			= 1 + \tfrac{(a-b)^{2}(3a(a-1)+3b(b-1)(K-1)+4abK)}{4(a+(K-1)b)^{2}(a(1-a)+b(1-b))}
			:= 1 + c_{a,b,K} \times \psi_{a,b,K},
		\end{equation}
		where $\psi_{a,b,K} := 3a(a-1)+3b(b-1)(K-1)+4abK$
		and $c_{a,b,K} > 0$.
		% $c_{a,b,K}:=\tfrac{(a-b)^{2}}{4(a+(K-1)b)^{2}(a(1-a)+b(1-b))}>0$
	\end{theorem}
	As in Table~\ref{table:K2_homo_relative_performance}, the function $\psi_{a,b,K}$ is the discriminating term that explicitly characterizes the relative performance of ASE and LSE.
	
	Here $\psi_{a,b,K}$ satisfies
	$(4ab-3(a-b^{2}))K < \psi_{a,b,K} < (4ab)K$, and there are explicit constants $c_{a,b}^{(1)}$ and $c_{a,b}^{(2)}$ depending only on $a$ and $b$ such that $\tfrac{1}{K}c_{a,b}^{(1)} < c_{a,b,K} \times \psi_{a,b,K} < \tfrac{1}{K}c_{a,b}^{(2)}$. Taking $a$ and $b$ to be fixed, these observations allow Eq.~(\ref{eq:rhoStar_K_homo_balanced_affinity}) to be summarized in terms of $K$ as
	\begin{equation}
		\rho^{\star} = 1 + \Theta_{a,b}(\tfrac{1}{K}),
	\end{equation}
	demonstrating that $\rho^{\star} \rightarrow 1$ as $K \rightarrow \infty$. In words, for the class of SBMs under consideration, ASE and LSE in a sense have asymptotically (in $K$) equivalent embedding performance (via $\rho^{\star}$). This amounts to a statement concerning a sequence of models with a necessarily growing number of vertices in order to ensure the underlying assumption of equal block sizes.
	
	Rewriting the level-set $\psi_{a,b,K}=0$, which holds if and only if $\rho^{\star} = 1$, yields the equation
	\begin{equation}
	\label{eq:convexCombParam}
		\left(\tfrac{1-a}{b}\right)\tfrac{1}{K} + \left(\tfrac{1-b}{a}\right)\tfrac{K-1}{K} = \tfrac{4}{3},
	\end{equation}
	together with the observation that ASE~$>$~LSE (resp.~ASE~$<$~LSE) when the left-hand side of Eq.~(\ref{eq:convexCombParam}) is less than (resp.~greater than) the value $\tfrac{4}{3}$. The above equation perhaps interestingly depicts a convex combination in $K$ of a reparameterization in terms of the variables $\tfrac{1-a}{b}$ and $\tfrac{1-b}{a}$, where the value $\tfrac{4}{3}$ is interpretable as a Chernoff-based information theoretic threshold.
	
	The observation that $\psi_{a,b,K} > (4ab-3(a-b^{2}))K$ in the context of Eq.~(\ref{eq:rhoStar_K_homo_balanced_affinity}) implies a sufficient condition for determining a parameter region in which ASE $>$ LSE \emph{uniformly in K}. Specifically, the condition $(4ab-3(a-b^{2})) > 0$, equivalently written as $\tfrac{a-b^{2}}{ab} < \tfrac{4}{3}$, ensures that $\psi_{a,b,K} > 0$ and hence that $\rho^{\star} > 1$.
	
	%RegionPlot[Evaluate[Table[(3 a (a - 1) + 3 b (b - 1) (k - 1) + 4 a b k <= 0 && b < a), {k, 2, 4}]], {a, 0, 1}, {b, 0, 1}, PlotPoints -> 100]

\begin{remark}[Detectability and phase transitions in random graph models]
	With respect to the random graph literature, the setting considered in this paper corresponds to a strong consistency regime (i.e.~exact recovery) in which the block membership of each individual vertex is recovered almost surely for graphs on $n$ vertices with $n \rightarrow \infty$.
	For different regimes where edge probabilities are allowed to decrease as a function of $n$, numerous deep and fascinating detectability and phase transition phenomena are known, some of which also employ Chernoff divergence and related considerations \citep{abbe2018survey}. In the context of homogeneous balanced affinity SBMs, the quantity $\textrm{SNR}:=\tfrac{(a-b)^{2}}{K(a+(K-1)b)}$ has been shown to function as an important information-theoretic signal-to-noise ratio. Here too the SNR appears, albeit with respect to $c_{a,b,K}$, in the sense that
	\begin{align*}
		c_{a,b,K}
		:=\tfrac{(a-b)^{2}}{4(a+(K-1)b)^{2}(a(1-a)+b(1-b))}
		&\equiv \left(\tfrac{(a-b)^{2}}{K(a+(K-1)b)}\right)\tilde{c}_{a,b,K}
	\end{align*}
	for some constant $\tilde{c}_{a,b,K}>0$. Perhaps more interestingly,
	\begin{align*}
		c_{a,b,K}
		&\equiv \tfrac{1}{4}\left(\tfrac{\lambda_{\textrm{min}}(\bB(K))}{\lambda_{\textrm{max}}(\bB(K))}\right)^2\left(\tfrac{1}{\sigma^{2}(\bB_{11}(K))+\sigma^{2}(\bB_{12}(K))}\right)
	\end{align*}
	where $\sigma^{2}(\bB_{ij}(K))$ is the edge variance corresponding to a pair of vertices in blocks $i$ and $j$, together with $\lambda_{\textrm{min}}(\bB(K))=a-b$ and $\lambda_{\textrm{max}}(\bB(K))=a+(K-1)b$, noting that the constant factor $\tfrac{1}{4}$ could just as easily be absorbed by redefining $\psi_{a,b,K}$. It may well prove fruitful to further investigate these observations in light of existing results.
	%Here, however, the SNR term does not arise when discriminating between the relative embedding performance of ASE and LSE.
\end{remark}

%%%%%%%%%%%%%%%%%%%%%%%%%%
%%%%%%%%%%%%%%%%%%%%%%%%%%
%%%%%%%%%%%%%%%%%%%%%%%%%%	

	\section{Discussion and Conclusions}
	\label{sec:Summary}
	Loosely speaking, Laplacian spectral embedding may be viewed as a degree-normalized version of adjacency spectral embedding in light of Eq.~(\ref{eq:Laplacian}). As such, our analysis complements existing literature that seeks to understand normalization in the context of spectral methods \citep{sarkar2015role,vonLuxburg2007tutorial}. Moreover, our work together with \cite{rubindelanchy2017generalized} addresses network models exhibiting indefinite geometry, an area that has received comparatively limited attention in the statistical network analysis literature. The ability of indefinite modeling considerations to reflect widely-observed disassortative community structure is encouraging and suggests future research activity in this and related directions.
	
	Core-periphery network structure, broadly construed, is demonstrably ubiquitous in real-world networks \citep{csermely2013structure,holme2005core,leskovec2009community}. With this understanding and the ability of the SBM to serve as a building block for hierarchically modeling complex network structure, our findings pertaining to spectral embedding for core-periphery structure may be of particular interest.
	
	This paper examines the information-theoretic relationship between the performance of two competing, widely-popular graph embeddings and subsequent vertex clustering with an eye towards underlying network model structure.
	The findings presented in Section~\ref{sec:ElucidatingNetworkStructure} support the claim that, for \emph{sparsity} interpreted as $\bB$ having entries that are small, loosely speaking, ``Laplacian spectral embedding favors relatively sparse graphs, whereas adjacency spectral embedding favors not-too-sparse graphs.'' Moreover, our results provide evidence in support of the claim that ``adjacency spectral embedding favors certain core-periphery network structure.'' Of course, caution must be exercised when making such general assertions, since the findings in this paper demonstrate intricate and nuanced functional relationships linking spectral embedding performance to network model structure. Nevertheless, we believe such summarative statements are both faithful and useful for conveying a high-level, macroscopic overview of the investigation presented in this work.
	
	%Our results extend existing, preliminary Chernoff-based spectral embedding analysis and are, to the best of our knowledge, among the first-ever results relating  the question of how the choice of embedding procedure both influences subsequent inference and couples with underlying network model structure.

%%%%%%%%%%%%%%%%%%%%%%%%%%
%%%%%%%%%%%%%%%%%%%%%%%%%%
%%%%%%%%%%%%%%%%%%%%%%%%%%
	\if0\blind
		{
		\section*{Acknowledgments}
		This work is partially supported by the XDATA and D3M programs of the Defense Advanced Research Projects Agency (DARPA) and by the Acheson J. Duncan Fund for the Advancement of Research in Statistics at Johns Hopkins University.
		Part of this work was done during visits by JC and CEP to the Isaac Newton Institute for Mathematical Sciences at the University of Cambridge (UK) under EPSCR grant EP/K032208/1. JC thanks Zachary Lubberts for productive discussions.
		} \fi
%%%%%%%%%%%%%%%%%%%%%%%%%%
%%%%%%%%%%%%%%%%%%%%%%%%%%
%%%%%%%%%%%%%%%%%%%%%%%%%%
	\newpage
	\section{Supplementary material}
	\label{sec:Appendix}
	\subsection{Latent position geometry}
	\label{sec:LatentPositionGeometry}
	All stochastic block models in Definition~\ref{def:SBM} can be formulated as instantiations of generalized random dot product graph models possessing inherent latent position (vector) structure. Earlier observations for the two-block SBM in Section~\ref{sec:ElucidatingNetworkStructure} are summarized in the following table, for which the implicit underlying vector $\bPi$ may be viewed as an additional parameter space dimension that weights the latent positions $\nu_{1}$ and $\nu_{2}$ by $\pi_{1}$ and $\pi_{2}$, respectively.
	
	\begin{center}
		\begin{tabular}{|l|l|}
			\hline
			Model geometry: & Canonical latent positions: \\ \hline 
			\hline
			Positive definite $\bB(a,b,c)$
			& $\nu_{1}=(\sqrt{a},0), \nu_{2}=(b/\sqrt{a},\sqrt{ac-b^{2}}/\sqrt{a})$ in $\R^{2}$\\
			Indefinite $\bB(a,b,c)$
			& $\nu_{1}=(\sqrt{a},0), \nu_{2}=(b/\sqrt{a},\sqrt{b^{2}-ac}/\sqrt{a})$ in $\R^{2}$ \\
			Rank one $\bB(p^{2}, pq, q^{2})$
			& $\nu_{1}=p, \nu_{2}=q$ in $\R$\\
			\hline
		\end{tabular}
	\end{center}
		
	For the homogeneous balanced affinity two-block network structure investigated in Section~\ref{sec:K2_homo_balanced}, the latent position geometry can be equivalently reparameterized as two vectors on the circle of radius $r := \sqrt{a}$ separated by the angle $\theta := \arccos(b/a)$.
	This behavior generalizes to the homogeneous balanced affinity $K$-block model.
	
	When $\bB \equiv \bB(K)\in(0,1)^{K \times K}$ has value $a$ on the main diagonal and value $b$ on the off-diagonal with $0<b<a<1$, we can write $\bB=\bX\bX^{\top}$ via the Cholesky decomposition, where $\bX$ has rows given by $\bX=[x_{1}|x_{2}|\dots|x_{K}]^{\top}$. For each $i\in[K]$ let the zero-dilation of the $\R^{K}$ vector $x_{i}$ be denoted by $x_{i}^{\circ}:=(x_{i},0)\in\R^{K+1}$. For $K=2,3,4$, $\bX$ is given by
	
	\begin{align}
	\bX(2) &:= \left[
	\begin{matrix}
	\sqrt{a} & 0 \\
	\frac{b}{\sqrt{a}} & \sqrt{\frac{(a-b)(a+b)}{a}}
	\end{matrix}
	\right] , \\
	\bX(3) &:= \left[
	\begin{matrix}
	\sqrt{a} & 0 & 0\\
	\frac{b}{\sqrt{a}} & \sqrt{\frac{(a-b)(a+b)}{a}} & 0\\
	\frac{b}{\sqrt{a}} & \sqrt{\frac{(a-b)(a+b)}{a}}\frac{b}{a+b} & \sqrt{\frac{(a-b)(a+2b)}{a+b}}
	\end{matrix}
	\right] , \label{eq:X2_X3_Cholesky}\\
	\bX(4) &:= \left[
	\begin{matrix}
	\sqrt{a} & 0 & 0 & 0\\
	\frac{b}{\sqrt{a}} & \sqrt{\frac{(a-b)(a+b)}{a}} & 0 & 0\\
	\frac{b}{\sqrt{a}} & \sqrt{\frac{(a-b)(a+b)}{a}}\frac{b}{a+b} & \sqrt{\frac{(a-b)(a+2b)}{a+b}} & 0\\
	\frac{b}{\sqrt{a}} & \sqrt{\frac{(a-b)(a+b)}{a}}\frac{b}{a+b} & \sqrt{\frac{(a-b)(a+2b)}{a+b}}\frac{b}{a+2b} & \sqrt{\frac{(a-b)(a+3b)}{a+2b}}
	\end{matrix}
	\right]. \label{eq:X4_Cholesky}
	\end{align}
	By induction, for $K \ge 3$, the entries of the vector $x_{K}$ are given by
	\begin{equation}
	\small
	x_{K} = \left(x_{K-1}^{1}, x_{K-1}^{2}, \dots, x_{K-1}^{K-2}, \left(\tfrac{b}{a+(K-2)b}\right)x_{K-1}^{K-1},\sqrt{\tfrac{(a-b)(a+(K-1)b)}{a+(K-2)b}}\right)^{\top} \in \R^{K}.
	\end{equation}
	
	Only $\bI_{0}^{K}$ and $\bI_{K-1}^{1}$ are necessary with respect to combining possible inner products on account of the sign-flip involving $a-b$. Beginning with the second row in each of the $\bX$ matrices, the first column of each matrix can be written in the more illuminating form $\sqrt{a}\tfrac{b}{a}$.	
	
	For this specific $K$-block model, symmetry with respect to equally-spaced vectors on the $\sqrt{a}$-radius sphere in $\R^{K}$ together with block membership balancedness translates into shared covariance structure such that Eq.~(\ref{eq:rhoStar}) reduces to Eq.~(\ref{eq:rhoStarSimple}). The first two rows of $\bX$ are ideal candidates to serve as canonical latent positions for subsequent computation, since these vectors are maximally sparse in the sense of having the fewest non-zero entries and merely become zero-inflated as a function of $K$. These geometric considerations are crucial in the subsequent proof of Theorem~\ref{thrm:Kblock_rhoStar}.
	
	%%%%%%%%%%%%%%%%%%%%%%%%%%
	%%%%%%%%%%%%%%%%%%%%%%%%%%
	%%%%%%%%%%%%%%%%%%%%%%%%%%

	\subsection{Analytic derivations for the two-block SBM}
	\label{sec:apdx:twoBlock}
	The value of $\rho^{\star}$ in Eq.~(\ref{eq:K2rhoStar}) for the homogeneous balanced two-block SBM can be computed by brute force; however, such an approach offers only limited insight and understanding of how the covariance structure in Theorem~\ref{thrm:GRDPG_CLT_ASE} and Theorem~\ref{thrm:GRDPG_CLT_LSE} interact to yield differences in relative spectral embedding performance as measured via Chernoff information. This section offers a different approach to understanding $\rho^{\star}$ as a covariance-based spectral quantity.
	
	The following lemma is a general matrix analysis observation that establishes a correspondence between the inverse of a convex combination of $2 \times 2$ matrices and the inverses of the original $2 \times 2$ matrices. The proof of Lemma~\ref{lem:twoBlockMatrixInversion} follows directly from elementary computations and is therefore omitted. Extending Lemma~\ref{lem:twoBlockMatrixInversion} to $n \times n$ invertible matrices is intractable in general.
	\begin{lemma}
	\label{lem:twoBlockMatrixInversion}
		Let $\bM_{0}, \bM_{1} \in \R^{2 \times 2}$ be two invertible matrices. For each $t\in[0,1]$ define the matrix $\bM_{t}:=(1-t)\bM_{0}+t\bM_{1}$. Provided $\bM_{t}$ is invertible, then the inverse matrix $\bM_{t}^{-1}$ can be expressed as
		\begin{equation}
		\label{eq:twoBlockMatrixInverse}
			\bM_{t}^{-1}
			\equiv \tfrac{(1-t)\bM_{0}^{-1}+\textnormal{det}(\bM_{1}\bM_{0}^{-1})t\bM_{1}^{-1}}{\textnormal{det}(\bM_{1}\bM_{0}^{-1})t^{2}+\textnormal{tr}(\bM_{1}\bM_{0}^{-1})t(1-t)+(1-t)^{2}}.
		\end{equation}
	\end{lemma}
	If, in the context of Lemma~\ref{lem:twoBlockMatrixInversion}, $\textrm{det}(\bM_{1}\bM_{0}^{-1})=1$, then Eq.~(\ref{eq:twoBlockMatrixInverse}) simplifies to
	\begin{equation*}
		\bM_{t}^{-1}
		\equiv \tfrac{(1-t)\bM_{0}^{-1}+t\bM_{1}^{-1}}{t^{2}+\textrm{tr}(\bM_{1}\bM_{0}^{-1})t(1-t)+(1-t)^{2}},
	\end{equation*}
	which is nearly a convex combination of the inverse matrices $\bM_{0}^{-1}$ and $\bM_{1}^{-1}$ modulo division by a degree two polynomial in the parameter $t$. If, in addition, $\textrm{tr}(\bM_{1}\bM_{0}^{-1})\neq -2$ (which always holds when $\bM_{0}$ and $\bM_{1}$ are both positive definite), then the inverse matrix at the value $t=\tfrac{1}{2}$ further simplifies to
	\begin{equation}
		\bM_{1/2}^{-1}
		\equiv \left(\tfrac{2}{2+\textrm{tr}(\bM_{1}\bM_{0}^{-1})}\right)\left(\bM_{0}^{-1}+\bM_{1}^{-1}\right).
	\end{equation}
	For the homogeneous balanced two-block SBM considered in Section~\ref{sec:K2_homo_balanced}, one can explicitly check that the above $\rmDet(\cdot)$ and $\rmTr(\cdot)$ conditions are satisfied. Moreover, the value $t^{\star}=\tfrac{1}{2}$ achieves the supremum in both the numerator and denominator of $\rho^{\star}$ in Eq.~(\ref{eq:rhoStar}). With these observations in hand, it follows by subsequent computations that for both the positive definite and indefinite regimes,
	
	\begin{align*}
		\rho^{\star}
		%\equiv \frac{\rho_{\tn{A}}^{\star}}{\rho_{\tn{L}}^{\star}}
		%= \frac{\underset{t\in(0,1)}{\tn{sup}}
		%	\left[	t(1-t)\|\nu_{1}-\nu_{2}\|_{\bSig_{1,2}^{-1}(t)}^{2}\right]}
		%{\underset{t\in(0,1)}{\tn{sup}}
		%	\left[t(1-t)\|\tilde{\nu}_{1}-\tilde{\nu}_{2}\|_{\tilde{\bSig}_{1,2}^{-1}(t)}^{2}
		%	\right]}
		&= \frac{\|\nu_{1}-\nu_{2}\|_{\bSig_{1,2}^{-1}(1/2)}^{2}}
		{\|\tilde{\nu}_{1}-\tilde{\nu}_{2}\|_{\tilde{\bSig}_{1,2}^{-1}(1/2)}^{2}}\\
		&= \left(\frac{(\tfrac{2}{2+\rmTr(\bSig(\nu_{1})\bSig^{-1}(\nu_{2}))})}{(\tfrac{2}{2+\rmTr(\tilde{\bSig}(\nu_{1})\tilde{\bSig}^{-1}(\nu_{2}))})}\right) \times 
		\left(\frac{(\nu_{1}-\nu_{2})^{\top}(\bSig^{-1}(\nu_{1})+\bSig^{-1}(\nu_{2}))(\nu_{1}-\nu_{2})}{(\tilde{\nu}_{1}-\tilde{\nu}_{2})^{\top}(\tilde{\bSig}^{-1}(\nu_{1})+\tilde{\bSig}^{-1}(\nu_{2}))(\tilde{\nu}_{1}-\tilde{\nu}_{2})}
		\right)\\
		&= \left(\frac{2+\rmTr(\tilde{\bSig}(\nu_{1})\tilde{\bSig}^{-1}(\nu_{2}))}{2+\rmTr(\bSig(\nu_{1})\bSig^{-1}(\nu_{2}))}\right) \times 1\\
		&= 1 + \frac{\rmTr(
			\tilde{\bSig}(\nu_{1})\tilde{\bSig}^{-1}(\nu_{2}))
			- \rmTr(\bSig(\nu_{1})\bSig^{-1}(\nu_{2}))}{2+\rmTr(\bSig(\nu_{1})\bSig^{-1}(\nu_{2}))}\\
		&= 1 + \frac{(a-b)^{2}(3a(a-1)+3b(b-1)+8ab)}{4(a+b)^{2}(a(1-a)+b(1-b))}. 		
	\end{align*}
	%Here we have exploited the simple nature of the scaled (resp.\ unscaled) latent positions $\nu_{i}$ (resp.\ $\tilde{\nu}_{i}$) as well as their straightforward relationship to one another.

	\subsection{Proof of Theorem~\ref{thrm:Kblock_rhoStar}} 
	This section is dedicated to proving Theorem~\ref{thrm:Kblock_rhoStar} for $K\ge2$ block SBMs exhibiting homogeneous balanced affinity structure. The proof is divided into two parts which separately evaluate the suprema in the numerator and denominator of $\rho^{\star}$ in Eq.~(\ref{eq:rhoStar}). By invoking underlying symmetries in latent space and the covariance structure of the ASE and LSE limit results, respectively, we shall leverage the (considerably simpler) ASE computations (numerator) when working with LSE (denominator). Simplifying the numerator and denominator yields the more easily interpretable (shifted) expression of $\rho^{\star}$ provided in Eq.~(\ref{eq:rhoStar_K_homo_balanced_affinity}).
	
		\emph{Proof:} First recall the discussion of latent space geometry in Section~\ref{sec:LatentPositionGeometry}, specifically that for the homogeneous balanced affinity $K$-block SBM, the canonical latent positions can be arranged row-wise as a lower-triangular matrix $\bX$ where each latent position vector has norm $\sqrt{a}$ and each pair of distinct latent position vectors has common inner-product $b$. This rotational symmetry implies rotational symmetry for the block-conditional covariance matrices in Theorems~\ref{thrm:GRDPG_CLT_ASE}--\ref{thrm:GRDPG_CLT_LSE}, and as such, the formulation of $\rho^{\star}$ in Eq.~(\ref{eq:rhoStar_K_homo_balanced_affinity}) can be reduced to simply working with the latent position pair $\{\nu_{1},\nu_{2}\}$ without loss of generality. This pair is attractive, since the non-zero entries of these vectors remain unchanged for all $K \ge 2$. One need only work with the standard inner product since $d^{-}=0$.
	
	\subsubsection{Proof of Theorem~\ref{thrm:Kblock_rhoStar}: ASE (numerator)}
	Let $g(x,X_{1}):=\langle x, X_{1} \rangle(1-\langle x, X_{1} \rangle)$
	and for $0 < t < 1$ define
	$g_{t}(x_{1},x_{2},X_{1}):= tg(x_{1},X_{1})+(1-t)g(x_{2},X_{1})$.
	By Theorem~\ref{thrm:GRDPG_CLT_ASE}, $\bSig(x)=\bDel^{-1}\mathbb{E}[g(x,X_{1})X_{1}X_{1}^{\top}]\bDel^{-1}$, and therefore
	$\bSig_{1,2}(t)
		:= t\bSig(\nu_{1})+(1-t)\bSig(\nu_{2})
		= \bDel^{-1}\mathbb{E}[g_{t}(\nu_{1},\nu_{2},X_{1})X_{1}X_{1}^{\top}]\bDel^{-1}$.
	Evaluating the inner expectation yields
	\begin{align*}
		\mathbb{E}[g_{t}(\nu_{1},\nu_{2},X_{1})X_{1}X_{1}^{\top}]
		&= \sum_{i=1}^{K}\tfrac{1}{K}(t\langle\nu_{1},\nu_{i}\rangle(1-\langle\nu_{1},\nu_{i}\rangle)+(1-t)\langle\nu_{2},\nu_{i}\rangle(1-\langle\nu_{2},\nu_{i}\rangle))\nu_{i}\nu_{i}^{\top}\\
		%&= \tfrac{1}{K}\left(ta(1-a)+(1-t)b(1-b)\right)\nu_{1}\nu_{1}^{\top}
		%+ \tfrac{1}{K}\left(tb(1-b)+(1-t)a(1-a)\right)\nu_{2}\nu_{2}^{\top}
		%+ \sum_{i=3}^{K}\tfrac{1}{K}b(1-b)\nu_{i}\nu_{i}^{\top}\\
		&= b(1-b)\bDel + \left(\tfrac{a(1-a)-b(1-b)}{K}\right)\left[t\nu_{1}\nu_{1}^{\top}+(1-t)\nu_{2}\nu_{2}^{\top}\right]\\
		&= b(1-b)\bDel + \bN(c_{0}\bD_{t})\bN^{\top},
	\end{align*}
	where $\bN := [\nu_{1}|\nu_{2}] \in \R^{K \times 2}$, $c_{0}:=\left(\tfrac{a(1-a)-b(1-b)}{K}\right)$, and $\bD_{t} := \textnormal{diag}(t,1-t)$.
	Clearly $c_{0}\bD_{t}$ is invertible, as is $\bDel$ since the underlying distribution $F$ is non-degenerate. Moreover, $\bX$ is also invertible since the $K$-block model under consideration is also rank $K$. The relation $\bX^{\top}\bX=K\bDel$ implies $\bDel^{-1}=K\bX^{-1}(\bX^{\top})^{-1}$ and therefore $\bX\bDel^{-1}\bX^{\top} = K\bI$, so $\nu_{i}^{\top}\bDel^{-1}\nu_{j} = K\mathbb{I}_{ij}$ where $\mathbb{I}_{ij}$ denotes the indicator function for indices $i$ and $j$. Thus,
	$(c_{0}\bD_{t})^{-1}+\frac{1}{b(1-b)}\bN^{\top}\bDel^{-1}\bN = (c_{0}\bD_{t})^{-1} + \frac{K}{b(1-b)}\bI$, which is also invertible. By an application of the Sherman--Morrison--Woodbury matrix inversion formula (\cite{horn2012matrix}, Section 0.7.4), then
	\begin{align*}
		\mathbb{E}[g_{t}(\nu_{1},\nu_{2},X_{1})X_{1}X_{1}^{\top}]^{-1}
		&= \left(b(1-b)\bDel + \bN(c_{0}\bD_{t})\bN^{\top}\right)^{-1}\\
		&= \left(\tfrac{1}{b(1-b)}\right)\bDel^{-1}
		- \left(\tfrac{1}{b(1-b)}\right)^{2}\bDel^{-1}\bN\left(\tfrac{1}{c_{0}}\bD_{t}^{-1}+\tfrac{K}{b(1-b)}\bI\right)^{-1}\bN^{\top}\bDel^{-1}.
	\end{align*}	
	For $\nu := \nu_{1}-\nu_{2} = {\scriptstyle \left(\tfrac{a-b}{\sqrt{a}}, -\sqrt{\tfrac{(a-b)(a+b)}{a}}, 0, \dots,0 \right)^{\top}}\in\R^{K}$, then $\nu^{\top}\bDel\nu=\frac{2}{K}(a-b)^{2}$ and $\bN^{\top}\nu=(a-b)(1,-1)^{\top} \in \R^{2}$. These observations together with subsequent computations yield the following chain of equalities.
	
	\begin{align*}
		\|\nu\|_{\bSig_{1,2}^{-1}(t)}^{2}
		&= \nu^{\top}\left(\bDel^{-1}\mathbb{E}[g_{t}(\nu_{1},\nu_{2},X_{1})X_{1}X_{1}^{\top}]\bDel^{-1}\right)^{-1}\nu\\
		&= \nu^{\top}(\bDel \mathbb{E}[g_{t}(\nu_{1},\nu_{2},X_{1})X_{1}X_{1}^{\top}]^{-1} \bDel)\nu\\
		&= \nu^{\top}\bDel\left(\tfrac{1}{b(1-b)}\bDel^{-1}
		- \left(\tfrac{1}{b(1-b)}\right)^{2}\bDel^{-1}\bN\left(\tfrac{1}{c_{0}}\bD_{t}^{-1}+\tfrac{K}{b(1-b)}\bI\right)^{-1}\bN^{\top}\bDel^{-1}\right)\bDel\nu\\
		&= \nu^{\top}\left(\tfrac{1}{b(1-b)}\bDel
		- \left(\tfrac{1}{b(1-b)}\right)^{2}\bN\left(\tfrac{1}{c_{0}}\bD_{t}^{-1}+\tfrac{K}{b(1-b)}\bI\right)^{-1}\bN^{\top}\right)\nu\\
		&= \left(\tfrac{1}{b(1-b)}\right)\nu^{\top}\bDel\nu
		- \left(\tfrac{1}{b(1-b)}\right)^{2}\nu^{\top}\bN\left(\tfrac{1}{c_{0}}\bD_{t}^{-1}+\tfrac{K}{b(1-b)}\bI\right)^{-1}\bN^{\top}\nu\\
		&= \left(\tfrac{2(a-b)^{2}}{b(1-b)K}\right) - \left(\tfrac{a-b}{b(1-b)}\right)^{2}(1,-1)\left(\tfrac{1}{c_{0}}\bD_{t}^{-1}+\tfrac{K}{b(1-b)}\bI\right)^{-1}(1,-1)^{\top}\\
		&= \left(\tfrac{2(a-b)^{2}}{b(1-b)K}\right) - \left(\tfrac{a-b}{b(1-b)}\right)^{2} \hspace{.2em} \rmTr\left(\left(\tfrac{1}{c_{0}}\bD_{t}^{-1}+\tfrac{K}{b(1-b)}\bI\right)^{-1}\right)\\
		&= \left(\tfrac{2(a-b)^{2}}{b(1-b)K}\right) - \left(\tfrac{a-b}{b(1-b)}\right)^{2}\left(\tfrac{(a(1-a)-b(1-b))b(1-b)t}{((a(1-a)-b(1-b))t+b(1-b))K} + \tfrac{(a(1-a)-b(1-b))b(1-b)(1-t)}{((a(1-a)-b(1-b))(1-t)+b(1-b))K} \right)\\
		&= \tfrac{(a-b)^{2}(a(a-1)+b(b-1))}{(a(1-a)+(a(a-1)-b(b-1))t)(b(b-1)+(a(a-1)-b(b-1))t)K}.
	\end{align*}
	In particular,
	\begin{equation}
	\label{eq:Kblock_ASE_numerator}
		\underset{t\in(0,1)}{\tn{sup}}\left[t(1-t)\|\nu\|_{\bSig_{1,2}^{-1}(t)}^{2}\right]
		= \tfrac{1}{K}\tfrac{(a-b)^{2}}{a(1-a)+b(1-b)},
	\end{equation}
	where by underlying symmetry the supremum is achieved at $t^{\star}=\tfrac{1}{2}$ over the entire parameter region $0<b<a<1$.

	\subsubsection{Proof of Theorem~\ref{thrm:Kblock_rhoStar}: LSE (denominator)}
	Recall that for this model $\indefI\equiv\bI_{d}$ since $d^{-}=0$. From Theorem~\ref{thrm:GRDPG_CLT_LSE} for LSE, the block conditional covariance matrix for each latent position $x$ can be written in the modified form
	\begin{equation*}
		\tilde{\bSig}(x)
		= \mathbb{E}\left[\left(\tfrac{g(x,X_{1})}{\langle x,\mu \rangle} \right) \left(\tfrac{\tilde{\bDel}^{-1}X_{1}}{\langle X_{1},\mu \rangle}-\tfrac{x}{2\langle x,\mu \rangle}\right)\left(\tfrac{\tilde{\bDel}^{-1}X_{1}}{\langle X_{1},\mu \rangle}-\tfrac{x}{2\langle x,\mu \rangle}\right)^{\top}\right].
	\end{equation*}
	
	We begin with several preliminary observations in order to define the quantities $c_{1}, c_{2}$, and $c_{3}$. Namely, for each latent position (row) $x$ of $\bX$,
	\begin{align}
		\langle x, \mu \rangle &= \left(\tfrac{a+(K-1)b}{K}\right) =: c_{1} \rlap{,} \label{eq:InnerXMu}\\
		\mathbb{E}[g(x,X_{1})] &= \left(\tfrac{a(1-a)+(K-1)b(1-b)}{K}\right) =: c_{2} \rlap{,} \label{eq:ExpectedG}\\
		\mathbb{E}[g(x,X_{1})X_{1}] &:= \left(\tfrac{a(1-a)-b(1-b)}{K}\right)x + b(1-b)\mu =: c_{3}x+b(1-b)\mu \rlap{.} \label{eq:ExpectedGX1}
	\end{align}
	Subsequent computations yield
	\begin{align*}
		\bDel x &= \left(\tfrac{a-b}{K}\right)x + b\mu \rlap{,} \\
		\left[\bDel-\left(\tfrac{a-b}{K}\right)\bI\right]xx^{\top} &= b \mu x^{\top} \rlap{,} \\
		\langle \bDel x,x \rangle &= \left(\tfrac{a^{2}+(K-1)b^{2}}{K}\right) \rlap{,} \\
		\tilde{\bDel} &\equiv \mathbb{E}\left[\tfrac{1}{\langle X_{1},\mu \rangle}X_{1}X_{1}^{\top}\right] = \tfrac{1}{c_{1}}\bDel \rlap{.}
	\end{align*}
	The above observations allow us to write $\tilde{\bSig}(x)$ as
	\begin{align*}
		&\mathbb{E}\left[\left(\tfrac{g(x,X_{1})}{\langle x,\mu \rangle} \right) \left(\tfrac{\tilde{\bDel}^{-1}X_{1}}{\langle X_{1},\mu \rangle}-\tfrac{x}{2\langle x,\mu \rangle}\right)\left(\tfrac{\tilde{\bDel}^{-1}X_{1}}{\langle X_{1},\mu \rangle}-\tfrac{x}{2\langle x,\mu \rangle}\right)^{\top}\right]\\
		&= \tilde{\bDel}^{-1}\mathbb{E}\left[\tfrac{g(x,X_{1})}{\langle x,\mu \rangle}\left(\tfrac{X_{1}}{\langle X_{1},\mu \rangle}-\tfrac{\tilde{\bDel}x}{2\langle x,\mu \rangle}\right)\left(\tfrac{X_{1}}{\langle X_{1},\mu \rangle}-\tfrac{\tilde{\bDel}x}{2\langle x,\mu \rangle}\right)^{\top}\right]\tilde{\bDel}^{-1} \\
		&= \tfrac{1}{c_{1}}\bDel^{-1}\mathbb{E}\left[g(x,X_{1})\left(X_{1}-\tfrac{1}{2c_{1}}\bDel x\right)\left(X_{1}-\tfrac{1}{2c_{1}}\bDel x\right)^{\top}\right]\bDel^{-1}.
	\end{align*}
	Expanding the term inside the expectation and applying linearity of expectation allows us to analyze each piece in turn. The first term in the expansion can be analyzed via the previous computations under ASE. For the second term,
	\begin{align*}
		\mathbb{E}\left[\tfrac{1}{2c_{1}}g(x,X_{1})X_{1}x^{\top}\bDel\right]
		&= \tfrac{1}{2c_{1}}\mathbb{E}[g(x,X_{1})X_{1}]x^{\top}\bDel\\
		&= \tfrac{1}{2c_{1}}\left(c_{3}xx^{\top}+b(1-b)\mu x^{\top}\right)\bDel\\
		&= \tfrac{1}{2c_{1}}\left(c_{3}xx^{\top}+(1-b)\left[\bDel-(\tfrac{a-b}{K})\bI\right]xx^{\top}\right)\bDel\\
		&= \left(\tfrac{1-b}{2c_{1}}\right) \bDel xx^{\top} \bDel + \left(\tfrac{Kc_{3}-(a-b)(1-b)}{2c_{1}K}\right)xx^{\top}\bDel\\
		&= \left(\tfrac{1-b}{2c_{1}}\right)\bDel xx^{\top} \bDel + \left(\tfrac{a(b-a)}{2c_{1}K}\right)xx^{\top}\bDel.
	\end{align*}
	Note that the transpose of this matrix corresponds to the third term in the implicit expansion of interest (not shown). Finally, the fourth term simply reduces to the form
	\begin{align*}
		\mathbb{E}\left[g(x,X_{1})\left(\tfrac{1}{2c_{1}}\bDel x\right)\left(\tfrac{1}{2c_{1}}\bDel x\right)^{\top}\right]
		&= c_{2}\left(\tfrac{1}{2c_{1}}\bDel x\right)\left(\tfrac{1}{2c_{1}}\bDel x\right)^{\top}
		= \left(\tfrac{c_{2}}{4c_{1}^{2}}\right) \bDel xx^{\top} \bDel.
	\end{align*}
	Thus,
	\begin{align*}
		& \mathbb{E}\left[g(x,X_{1})\left(X_{1}-\tfrac{1}{2c_{1}}\bDel x\right)\left(X_{1}-\tfrac{1}{2c_{1}}\bDel x\right)^{\top}\right]\\
		&= \mathbb{E}[g(x,X_{1})X_{1}X_{1}^{\top}]
		- \mathbb{E}[\tfrac{1}{2c_{1}}g(x,X_{1})X_{1}x^{\top}\bDel]
		- \mathbb{E}[\tfrac{1}{2c_{1}}g(x,X_{1})X_{1}x^{\top}\bDel]^{\top}
		+ \mathbb{E}[g(x,X_{1})(\tfrac{1}{2c_{1}}\bDel x)(\tfrac{1}{2c_{1}}\bDel x)^{\top}]\\
		&= \mathbb{E}[g(x,X_{1})X_{1}X_{1}^{\top}]
		- \left(\tfrac{a(b-a)}{2c_{1}K}\right)xx^{\top}\bDel 
		- \left(\tfrac{a(b-a)}{2c_{1}K}\right)\bDel xx^{\top}
		+ \left(\tfrac{c_{2}}{4c_{1}^{2}}-\tfrac{1-b}{c_{1}}\right)\bDel xx^{\top} \bDel.
		%&= \mathbb{E}[g(x,X_{1})X_{1}X_{1}^{\top}]
		%+ \tfrac{a(a-b)}{2c_{1}K}xx^{\top}\bDel 
		%+ \tfrac{a(a-b)}{2c_{1}K}\bDel xx^{\top}
		%+ \left(\tfrac{c_{2}-4c_{1}(1-b)}{4c_{1}^{2}}\right)\bDel xx^{\top} \bDel.
	\end{align*}
	Let $\bM_{1}\equiv\bM_{1}(t):=\bN\bD_{t}\bN^{\top}$ and  $\bM_{2}:=\bDel$ with respect to the notation introduced earlier in the derivation for ASE. By completing the appropriate matrix product, there are explicit constants $\{d_{i}\}_{i=1}^{4}$ depending on $a$, $b$, and $K$, such that
	\begin{align*}
		\tilde{\bSig}_{1,2}(t)
		&= t\tilde{\bSig}(\nu_{1})+(1-t)\tilde{\bSig}(\nu_{2})\\
		&= \bDel^{-1}\left(d_{1}\bDel+d_{2}\bN\bD_{t}\bN^{\top} + d_{3}\bN\bD_{t}\bN^{\top}\bDel + d_{3}\bDel\bN\bD_{t}\bN^{\top} + d_{4}\bDel\bN\bD_{t}\bN^{\top}\bDel \right)\bDel^{-1}\\
		&= \bDel^{-1}\left(\left[d_{1}\bM_{2}+d_{5}\bM_{1}\right] + (\bI+d_{6}\bM_{2})(d_{7}\bM_{1})(\bI+d_{6}\bM_{2})\right)\bDel^{-1}\\
		&=: \bDel^{-1}\left(\bM_{3} + \bM_{4}\right)\bDel^{-1},
	\end{align*}
	where $\bM_{3}\equiv\bM_{3}(t):=d_{1}\bM_{2}+d_{5}\bM_{1}(t)$ and $\bM_{4}\equiv\bM_{4}(t)=(\bI+d_{6}\bM_{2})(d_{7}\bM_{1}(t))(\bI+d_{6}\bM_{2})$.
	
	Note that
	$\tilde{\nu}_{k}
	:=\left(\tfrac{1}{\langle \nu_{k}, \mu \rangle}\right)^{1/2} \times \nu_{k}
	= \left(\tfrac{K}{a+(K-1)b}\right)^{1/2} \times \nu_{k}$ for $k=1,2$, so
	\begin{align*}
		\|\tilde{\nu}\|_{\tilde{\bSig}_{1,2}^{-1}(t)}^{2}
		&= \tilde{\nu}^{\top}\tilde{\bSig}_{1,2}^{-1}(t)\tilde{\nu}
		= \left(\tfrac{K}{a+(K-1)b}\right)\nu^{\top}\bDel\left(\bM_{3}+\bM_{4}\right)^{-1}\bDel\nu.
	\end{align*}
	The above matrix inversion can again be carried out via the Sherman--Morrison--Woodbury formula. We omit the algebraic details. Subsequent computations and simplification yield
	\begin{equation}
	\label{eq:Kblock_LSE_denomoinator}
		\underset{t\in(0,1)}{\tn{sup}}\left[t(1-t)\|\tilde{\nu}\|_{\tilde{\bSig}_{1,2}^{-1}(t)}^{2}\right]
		= \tfrac{4(a-b)^{2}(a+(K-1)b)^{2}}{4(a(1-a)+b(1-b))(a+(K-1)b)^{2}K+(a-b)^{2}K(3a(a-1)+3b(b-1)(K-1)+4abK)}
	\end{equation}
	where by underlying symmetry the supremum is achieved at $t^{\star}=\tfrac{1}{2}$ over the entire parameter region $0<b<a<1$. Taken together, Eq.~(\ref{eq:Kblock_ASE_numerator}) and Eq.~(\ref{eq:Kblock_LSE_denomoinator}) simplify to yield $\rho^{\star}$ as in Eq.~(\ref{eq:rhoStar_K_homo_balanced_affinity}), thereby completing the proof. $\square$

\newpage
\bibliographystyle{chicago}
\bibliography{main}

\end{document}